\documentclass[11pt]{article}
\usepackage{afterpage}

\usepackage[utf8]{inputenc}
\usepackage{amsfonts,amssymb,amsmath,amsthm,latexsym,epsfig,amscd,bbm,stmaryrd,mathrsfs}
\usepackage[english]{babel}
\usepackage{graphicx}
\usepackage{color}
\usepackage{tikz}
\usepackage{courier}
\usepackage{bbm}
\usepackage{enumerate}
\usepackage{psfrag}
\usepackage{wasysym}
\usepackage{type1cm}
\usepackage{dsfont}
\usepackage{geometry}
\usepackage{marginnote}

\numberwithin{equation}{section}

\theoremstyle{plain} 

\newtheorem{theorem}{Theorem}[section]
\newtheorem{corollary}[theorem]{Corollary}
\newtheorem{lemma}[theorem]{Lemma}
\newtheorem{proposition}[theorem]{Proposition}
\newtheorem{definition}[theorem]{Definition}
\newtheorem{assumption}[theorem]{Assumption}

\newcommand{\ssup}[1] {{{\scriptscriptstyle{({#1}})}}} 
\newcommand{\EE}{\mathbb{E}}

\newcommand{\N}{\mathbb{N}}

\newcommand{\R}{\mathbb{R}}

\newcommand{\cC}{\mathcal{C}}
\newcommand{\cB}{\mathcal{B}}

\newcommand{\cO}{\mathcal{O}} 
\newcommand{\cT}{\mathcal{T}} 
 
\newcommand{\Prob}{\mathrm{P}}

\newcommand{\eee}{\mathrm{e}}

\def \cL {{\mathcal L}}

\DeclareMathOperator\supp{supp}
\DeclareMathOperator\dist{dist}
\DeclareMathOperator\diam{diam}

\newcommand{\scrP}{\mathscr{P}}

\newcommand{\BB}{\mathcal B}
\newcommand{\CC}{\mathcal C}

\newcommand{\MM}{\mathcal M}
\newcommand{\NN}{\mathcal N}

\newcommand{\YY}{\mathcal Y}

\def\N{\mathbb{N}}
\def\R{\mathbb{R}}

\def\P{\mathbb{P}}
\def\E{\mathbb{E}}

\newcommand{\cc}{{\text{\rm c}}}
\newcommand{\texte}{\text{\rm e}}

\def \GW {{\mathcal{G}\mathcal{W}}}

\newcommand{\mr}{\mathring}
\newcommand{\textd}{\text{\rm d}\mkern0.5mu}

\def\1{{\mathchoice {1\mskip-4mu\mathrm l}
{1\mskip-4mu\mathrm l} 
{1\mskip-4.5mu\mathrm l} {1\mskip-5mu\mathrm l}}}

\def \ba {\begin{array}}
\def \ea {\end{array}}

\def \P  {{\mathbb P}}
\def \E  {{\mathbb E}}

\def \cL {{\mathcal L}}

\def \cN {{\mathcal N}}
\def \cO {{\mathcal O}}
\def \cP {{\mathcal P}}

\def \cT {{\mathcal T}}
\def \GW {{\mathcal{G}\mathcal{W}}}

\def \e {\mathrm{e}}
\def \ee {\mathrm{e}}
\def \dd {\mathrm{d}}

\DeclareSymbolFont{symbolsC}{U}{pxsyc}{m}{n}
\DeclareMathSymbol{\opentimes}{\mathrel}{symbolsC}{93}

\newcommand{\Probgr}{\mathfrak{P}}

\title{The parabolic Anderson model on a Galton-Watson tree with normalised Laplacian}

\author{

D. Wang
\footnotemark[1]
}

\footnotetext[1]{
Mathematical Institute, Leiden University, P.O.\ Box 9512, 2300 RA Leiden, The Netherlands.
\\
Email: {\tt d.wang@math.leidenuniuv.nl}}

\date{\today}

\begin{document}

\maketitle

\begin{abstract}
In earlier work by den Hollander, K\"onig, and dos Santos, the asymptotics of the total mass of the solution to the parabolic Anderson model was studied on an almost surely infinite Galton-Watson tree with an i.i.d. potential having a double-exponential distribution. The second-order contribution to this asymptotics was identified in terms of a variational formula that gives information about the local structure of the region where the solution is concentrated. 

The present paper extends this work to the degree-normalised Laplacian. The normalisation causes the Laplacian to be non-symmetric and which leads to different spectral properties. We find that the leading order asymptotics of the total mass remains the same, while the second-order correction coming from the variational formula is different. We also find that the optimiser of the variational formula is again an infinite tree with minimal degrees. Both of these results are shown to hold under much milder conditions than for the regular Laplacian. 

\vspace{0.5cm}



\medskip\noindent
{\bf Acknowledgments:}
The work in this paper was supported through Gravitation-grant NETWORKS-024.002.003 of the Netherlands Organisation for Scientific Research (NWO).\\
The author would like to thank Frank den Hollander for suggesting the PAM with the normalised Laplacian and for helpful discussions. The author would also like to thank Michal Bassan for helpful discussions.
\normalsize
\end{abstract}

\newpage


\newpage


\section{Introduction and main results}
\label{s:intro}

\subsection{The PAM and intermittency}
\label{ss:introPAM}

\noindent The \emph{parabolic Anderson model} (PAM) on a graph $G = (V,E)$, is the Cauchy problem for the heat equation with a random potential:
\begin{equation}
\label{e:PAMdef}
\begin{array}{llll}
\partial_t u(x,t) &=& (\Delta u)(x,t) + \xi(x) u(x,t), &x \in V,\, t>0,\\
u(x,0) &=& \delta_\cO(x), &x \in V,
\end{array}
\end{equation} 
where $\cO$ is a vertex in $V$, $\{\xi(x)\}_{x\in V}$ is the random potential defined on $V$, and $\Delta$ is the normalised discrete Laplacian, defined by
\begin{equation}
\label{e:Laplacian}
(\Delta f)(x) := \frac{1}{\deg(x)}\sum_{{y\in V:} \atop { \{x,y\} \in E} } [f(y) - f(x)], \qquad x \in V,\,f\colon\,V\to\R.
\end{equation}

Most of the literature has focused on $\mathbb{Z}^d$ and we refer the reader to \cite{K2016} for a comprehensive study. Other choices include the complete graph \cite{FM1990}, the hypercube \cite{AGH2020}, and more recently, the regular tree \cite{dhW2023} and the Galton-Watson tree \cite{dHKdS2020}, \cite{dhW2022} and \cite{AP2023}. The PAM can also be studied on continuous spaces, the most extensively studied being $\mathbb{R}^d$. On such spaces, the Laplacian and the potential are defined analogously and we again refer the reader to \cite{K2016} for more background.

The PAM may be interpreted as a system of particles such that particles are killed with rate $\xi^-(x)$ or are split into two with rate $\xi^+(x)$ at every vertex $x$. At the same time, each particle jumps independently with $\Delta$ as generator. The solution $u(x,t)$ can be interpreted as the number of particles or \textit{mass} present at vertex $x$ at time $t$ when the initial condition at time $0$ is $\delta_\cO(x)$. See \cite[Section 1.2]{GM1990} for further details.

The solution is known to exhibit a phenomenon known as \textit{intermittency}, meaning that the solution concentrates on small regions of the graph known as \textit{intermittent islands}. This is well studied on $\mathbb{Z}^d$, where it is known that the sizes of the island(s) depend on the distribution of the tail of the potential, and can be separated into four classes (see \cite[Section 3.4]{K2016}). Throughout the paper, we work in the \textit{double-exponential} class where the potential $\xi = \{\xi(x)\}_{x\in V}$ consists of i.i.d. random variables satisfying
\begin{assumption}{\bf [Asymptotic double-exponential potential]} 
\label{ass:pot}
$\mbox{}$\\
For some $\varrho \in (0,\infty)$,
\begin{equation}
\label{e:DE}
\Prob \left( \xi(0) \geq 0\right) = 1, \qquad \Prob \left(\xi(0) > u \right) = \ee^{-\ee^{u/\varrho}} \;\; 
\text{for } u \text{ large enough.}
\end{equation}
\end{assumption}
\noindent The main feature of this choice is that the intermittent islands are not single vertices, whilst also having the property that their sizes do not change with time - a critical fact in our analysis.

\subsection{The PAM on a Galton-Watson tree}
\label{ss:GW}
The present paper analyses the PAM on the graph generated by the Galton-Watson process. The graph is generated by taking a root $\cO$ and attaching $D$ vertices (known as offspring), where $D$ is a random variable. Each offspring has $D-1$ offspring attached to it, where $D$ is an identically distributed but independent copy of $D$. This is repeated forever, or until the process dies out. Let $\GW = (V,E,\cO)$ be the resulting graph and let $\mathfrak{P}$ and $\mathfrak{E}$ denote probability and expectation with respect to $\GW$. Similarly, let $\mathcal{P}$ and $\mathcal{E}$ denote probability and expectation with respect to $D$.
\begin{assumption}{\bf [Exponential tails]}
\label{ass:deg}
$\mbox{}$\\
{\rm (1)} $d_{\min} := \min \supp(D) \geq 2$ and $\mathcal{E}[D] \in (2,\infty)$.\\
{\rm (2)} $\mathcal{E}\big[\ee^{aD}\big] < \infty$ for all $a \in (0,\infty)$.
\end{assumption}

\noindent
Under this assumption, $\GW$ is $\Probgr$-a.s.\ an infinite tree. Moreover,
\begin{equation}
\label{e:volumerateGW}
\lim_{r \to \infty} \frac{\log |B_r(\cO)|}{r} = \log \mathcal{E}[D] =: \vartheta \in (0,\infty) \qquad \Probgr-a.s.,
\end{equation}
where $B_r(\cO) \subset V$ is the ball of radius $r$ around $\cO$ in the graph distance (see e.g.\ \cite[pp.~134--135]{LP2016}). Note that this ball depends on $\GW$ and is therefore random.\\

The proper choice of Laplacian depends on the setting. In the case of the complete graph and the hypercube, when the limit of the number of vertices going to infinity is taken, only the normalised Laplacian gives a meaningful limit. In the case of regular graphs such as $\mathbb{Z}^d$ and the regular tree, normalising the Laplacian by the degree simply amounts to rescaling time, and both the techniques and the results can be easily inferred accordingly. We will focus on the Galton-Watson tree, which is not only inhomogeneous but is also random, and hence the choice of Laplacian does play a role. The present paper considers the normalised Laplacian, in contrast to \cite{dHKdS2020}, \cite{dhW2022}, and \cite{AP2023}, and investigate how this choice affects the results and methods used in $\cite{dHKdS2020}$.

Under Assumptions~\ref{ass:pot} and \ref{ass:deg}, the criteria for existence and uniqueness of a non-negative solution of \eqref{e:PAMdef} are met (see \cite{GM1990} and \cite[Appendix C]{dhW2022}) and is given by the well-known Feynman-Kac representation. With the choice of initial condition in \eqref{e:PAMdef} and Laplacian in \eqref{e:Laplacian}, this amounts to
\begin{equation}
\label{e:FK}
u(x,t) = \E_\cO \left[ \exp \left\{\int_0^t \xi(X_s) \textd s \right\}\,\1\{X_t = x\} \right],
\end{equation}
where $X=(X_t)_{t \geq 0}$ is the continuous-time random walk on the vertices $V$ with jump rate $1$ on each vertex (or equivalently with jump rate equal to the inverse of the degree along the edges $E$), and $\P_\cO$ denotes the law of $X$ given $X_0=\cO$. The quantity we will be interested in is the \textit{total mass}, given by
\begin{equation}
\label{e:mass} 
U(t):= \sum_{x\in V} u(x,t) = \E_\cO \left[ \exp \left\{\int_0^t \xi(X_s) \textd s \right\}\right],
\end{equation}
in particular, its asymptotics as $t \to \infty$. 

An important distinction is made between the \textit{quenched} and the \textit{annealed} total mass, i.e. the total mass taken almost surely with respect to or averaged over the sources of randomness, respectively. Since this paper aims to follow the framework of \cite{dHKdS2020}, we also consider the quenched setting for both the graph and potential. We refer to \cite{dhW2023} for corresponding results for the annealed total mass on a regular tree. The annealed setting for the Galton-Watson tree (averaged over just the potential or over both the graph and potential) remain open.

\subsection{Main results and discussion}
\label{ss:resdis}
To state our results, we first introduce some quantities of interest pertaining to the \textit{characteristic variational formula} associated with Assumption~\ref{ass:pot}. The latter describes the shape and profile of the solution in the intermittent islands and captures the second-order asymptotics of the total mass. We refer to \cite{K2016} for more details on the variational formula and its relationship with the PAM.

Denote the set of probability measures on $V$ by $\cP(V)$. For $p \in \cP(V)$, define
\begin{equation}
\label{e:defIJ}
I_E(p) := \sum_{\{x,y\} \in E} \left(\sqrt{\tfrac{p(x)}{\deg(x)}}-\sqrt{\tfrac{p(y)}{\deg(y)}}\right)^2,
\qquad J_V(p) := - \sum_{x \in V} p(x) \log p(x),
\end{equation}
and set
\begin{equation}
\label{e:defchiG}
\chi_G(\varrho) := \inf_{p \in \cP(V)} [I_E(p) + \varrho J_V(p)], \qquad \varrho \in (0,\infty).
\end{equation}
As shown later, the first term arises from the Laplacian and coincides with the large deviation rate function for the empirical distribution of the random walk in \eqref{e:FK}, while the second term comes from the choice of the double-exponential potential. Furthermore, define the constant
\begin{equation}
\label{e:deftildechi}
\widetilde{\chi}(\varrho) := \inf \big\{ \chi_T(\varrho) \colon\, T \text{ is an infinite tree with degrees in } \supp(D) \big\},
\end{equation}
with $\chi_G(\varrho)$ defined in \eqref{e:defchiG}, and abbreviate
\begin{equation}
\label{mathfrakrdef}
\mathfrak{r}_t = \frac{\varrho t}{\log\log t}.
\end{equation}
\begin{theorem}{\bf [Total mass asymptotics]}
\label{t:QLyapGWT}
Subject to Assumptions~{\rm \ref{ass:pot}}--{\rm \ref{ass:deg}},
\begin{equation}
\label{e:QLyapGWT}
\frac{1}{t} \log U(t) = \varrho \log (\vartheta\mathfrak{r}_t) 
-\varrho - \widetilde{\chi}(\varrho) + o(1), \quad t \to \infty, 
\qquad (\Prob \times \Probgr)\text{-a.s.}
\end{equation}
\end{theorem}
\noindent The proof of Theorem~\ref{t:QLyapGWT} is given in Section~\ref{s.proofmain}. For a heuristic explanation on how the terms in \eqref{e:QLyapGWT} arise and how they relate to the asymptotics of the total mass, we refer the reader to \cite[Section 1.5]{dHKdS2020}.

For $d \ge 2$, let $\cT_d$ denote the \emph{infinite homogeneous tree} with degree equal to $d$ at every vertex.
\begin{theorem}{\bf [Identification of the minimiser]}
\label{t:tildechilargerho}
If $\varrho \geq \tfrac{1}{(d_{\rm min}-1)\log(d_{\rm min}+1)}$, then $\widetilde{\chi}(\varrho) = \chi_{\cT_{d_{\min}}}(\varrho)$.
\end{theorem}
\noindent The proof of Theorem~\ref{t:tildechilargerho} is given in Section \ref{ss:theoremvarproof}.

Comparing our results to those obtained in \cite{dHKdS2020} and \cite{dhW2022}, we see that the choice of Laplacian indeed has an effect, albeit in a subtle way. The leading-order terms in \eqref{e:QLyapGWT} remain unchanged, while the second-order (variational formula) term stemming from \eqref{e:defchiG} is different, due to $I_E$ in \eqref{e:defIJ} being normalised by the degrees. This normalisation was not present in \cite{dHKdS2020} and \cite{dhW2022}. In addition, normalising the Laplacian results in a `slow down' of the random walk in \eqref{e:FK} compared with the analogous formula in \cite{dHKdS2020}. As will be shown later on, this leads to simplifications in several key lemmas and leads to Theorem~\ref{t:QLyapGWT} holding under the milder tail condition in Assumption~\ref{ass:deg}(2).

The different Laplacian and $I_E$ function have surprisingly minimal effects on Theorem~\ref{t:tildechilargerho}: the optimal tree is still $\cT_{d_{\min}}$, exactly as was found in \cite{dHKdS2020}. The main difference is that our result holds for a greater range of $\varrho$ values compared to \cite{dHKdS2020}, which required the sharper restriction $\varrho \geq 1/\log(d_{\rm min}+1)$. We believe that the minimal tree is the minimiser for all $\varrho$ and that it is also the unique minimiser, however this remains open. It is also worth noting that the object $\widetilde{\chi}(\varrho)$ is well understood. The case $d_{\min}= 2$ corresponds to the the variational problem on $\mathbb{Z}$ and has been studied in \cite{dHG1999}. For $d_{\min} > 2$ we refer the reader to \cite{dhW2023}, where the variational formula 
\begin{equation*}
\bar{\chi}_{d_{\min}}(\varrho) =  \inf_{p \in \cP(V)} [\bar{I}_E(p) + \varrho J_V(p)], \quad  \bar{I}_E(p) = \sum_{\{x,y\} \in E} \left( \sqrt{p(x)} - \sqrt{p(y)}\,\right)^2,
\end{equation*}
was studied. Clearly,
\begin{equation*}
\widetilde{\chi}(\varrho) = \tfrac{1}{d_{\min}}\bar{\chi}_{d_{\min}}(d_{\min}\varrho).    
\end{equation*}
\paragraph{Outline.} The remainder of the paper is dedicated to the proof of Theorems~\ref{t:QLyapGWT} and \ref{t:tildechilargerho}, and follows the framework developed in \cite{dHKdS2020}. Section 2.1 is novel and deals with the spectral estimates of the Anderson Hamiltonian $\Delta + \xi$, which are different due to $\Delta$ no longer being symmetric with respect to the usual inner product. Sections 2.2 and 2.3 collect the necessary results regarding the Galton-Watson tree and the potential from \cite{dHKdS2020} and \cite{dhW2022}. All of these results carry over directly since the Laplacian and random walk play no role. Section 3 follows the \textit{path expansion} technique from \cite{dHKdS2020} and adapts the results to the random walk in \eqref{e:FK}. Section 4 is dedicated to the proof of Theorem~\ref{t:QLyapGWT}, and follows \cite{dHKdS2020}. Section 5 deals with the analysis of the variational formula \eqref{e:deftildechi} including the proof of Theorem~\ref{t:tildechilargerho}, which applies the gluing argument from \cite{dHKdS2020}.

\section{Preliminaries}
In this section we collect results that are needed later. Section~\ref{prespectral} investigates how the normalisation affects the spectral properties of the Laplacian. Section~\ref{s:struct} collects two vital facts about the Galton-Watson tree. Section~\ref{ss:potestim} collects results regarding the potential. 
\subsection{Related Spectral Problems}
\label{prespectral}
We introduce an alternative representation for $\chi$ in \eqref{e:defchiG} in terms of a \lq dual\rq\ variational formula. Fix $\varrho \in (0,\infty)$ and a graph $G=(V,E)$. The functional
\begin{equation}
\label{e:cL}
\cL(q;G) := \sum_{x \in V} \ee^{q(x)/\varrho}\in[0,\infty], \qquad q\colon\, V \to [-\infty,\infty), 
\end{equation}
plays the role of a large deviation rate function for the potential $\xi$ in $V$ (compare with \eqref{e:DE}). For $\Lambda \subset V$, define
\begin{equation}
\label{e:defhatchi}
\widehat{\chi}_{\Lambda}(G) := - \sup_{\substack{q\colon V \to [-\infty,\infty), \\ \cL(q;G) \leq 1}} 
\lambda_{\Lambda}(q;G)\in[0,\infty),
\end{equation}
where $\lambda_{\Lambda}(q;G)$ is the principal eigenvalue of the Anderson Hamiltonian $\Delta + q$ on the set $\Lambda$ with zero boundary condition. The condition $\cL(q;G) \leq 1$ under the supremum ensures that the potentials $q$ have a fair probability under the i.i.d.\ double-exponential distribution. \\

\begin{proposition}{\bf [Alternative representations for $\chi$]}
\label{p:dualrepchi}
For any graph $G = (V,E)$ and any $\Lambda \subset V$,\\
\begin{equation}
\widehat{\chi}_\Lambda(\varrho;G) \geq \widehat{\chi}_V(\varrho; G)= \chi_G(\varrho).
\end{equation}
\end{proposition}
\noindent Proposition~\ref{p:dualrepchi} is not essential in the proof of Theorem~\ref{t:QLyapGWT}, but is stated here to provide additional context to some of the results below. The proof is given in Section~\ref{proofsection}.\\

Recall the Rayleigh-Ritz formula for the principal eigenvalue $\lambda_\Lambda(q;G)$,
\begin{equation}
\label{e:RRformula}
\begin{aligned}
\lambda_\Lambda(q;G) = \sup \big\{ \langle (\Delta + q) \phi, \phi \rangle 
\colon\, \phi \in \R^{V}, \,\supp \phi \subset \Lambda, \, \|\phi\|=1 \big\}.
\end{aligned}
\end{equation}
As alluded to in Section~\ref{ss:resdis}, $\Delta$ (and therefore also $\Delta+q$) is not symmetric with respect to the usual $\ell^2$ inner product, but is symmetric with respect to the degree-weighted inner product
\begin{equation}
\label{e:definnerprod}
\langle \phi,\psi \rangle:= \sum_{x\in \Lambda} \deg(x) \phi(x)\psi(x),
\end{equation}
in the sense that $\langle \Delta \phi, \psi \rangle = \langle \phi, \Delta \psi\rangle$ (see \cite[Section 2]{JJ2002} for further details). Henceforward, all inner products will be with respect to \eqref{e:definnerprod}.
\begin{lemma}{\bf [Spectral bounds]}
\label{lem:specbd}
\begin{enumerate}
\item[{\rm (1)}] 
For any $\Gamma \subset \Lambda \subset V$, 
\begin{equation}
\label{e:monot_princev}
\max_{z \in \Gamma} q(z) - 1 \leq \lambda_\Gamma(q;G) 
\leq \lambda_\Lambda(q;G) \leq \max_{z \in \Lambda} q(z).
\end{equation}
\item[{\rm (2)}] 
The eigenfunction corresponding to $\lambda_\Lambda(q;G)$ can be taken to be non-negative.
\item[{\rm (3)}] 
If $q$ is real-valued and $\Gamma \subsetneq \Lambda$ is finite and connected in $G$, then the second inequality in \eqref{e:monot_princev} is strict and the eigenfunction corresponding to $\lambda_\Lambda(q;G)$ is strictly positive.
\end{enumerate}
\end{lemma}
\begin{proof}
Write
\begin{equation}
\begin{aligned}
\langle (\Delta + q) \phi, \phi \rangle = -\sum_{\{x, y\}\in E_{\Lambda}}[\phi(x)-\phi(y)]^{2} + \sum_{x\in \Lambda} 
 \deg(x) q(x)\phi(x)^{2}.   
\end{aligned}
\end{equation}
The upper bound in \eqref{e:monot_princev} follows from the estimate
\begin{align*}
\langle (\Delta + q) \phi, \phi \rangle & \leq \sum_{x \in \Lambda} \deg(x)q(x)\phi(x)^2 \leq \max_{z \in \Lambda} q(z)\sum_{x \in \Lambda}\deg(x)\phi(x)^2  = \max_{z \in \Lambda} q(z).
\end{align*}
To get the lower bound in \eqref{e:monot_princev}, we use the fact that $\lambda_\Lambda$ is non-decreasing in $q$. Let $\Bar{z} = \arg \max q(z)$.
Replacing $q(z)$ by $-\infty$ for every $z \neq \bar z$ and taking the test function $\bar\phi = \tfrac{1}{\sqrt{\deg(\bar{z})}}\delta_{\bar z}$, we get that 
\begin{equation}
\begin{aligned}
\lambda_\Lambda(q;G) 
& \geq - \sum_{ {x, y\in \Lambda:} \atop { \{x,y\} \in E_\Lambda} }\left[\bar \phi(x)-\bar \phi(y)\right]^2
+ \sum_{x \in \Lambda}\deg(x) q(x)\bar \phi(x)^2\\
&=-\sum_{ {y \in \Lambda:} \atop { \{\bar z,y\} \in E_\Lambda} } \frac{1}{\deg(\bar{z})} + q(\bar z)
= -1 + \max_{z \in \Lambda} q(z),  
\end{aligned}
\end{equation}
which settles the claim in (1). The claims in (2) and (3) are standard.
\end{proof}

\noindent Inside $\GW$, fix a finite connected subset $\Lambda \subset V$, and let $H_\Lambda$ denote the Anderson Hamiltonian in $\Lambda$ with zero Dirichlet boundary conditions on $\Lambda^c = V \backslash \Lambda$ (i.e. the restriction of the operator $H_G = \Delta + \xi$ to the class of functions supported on $\Lambda$). For $y \in \Lambda$, let $u^y_\Lambda$ be the solution of
\begin{equation}
\label{e:PAMalt}
\begin{array}{llll}
\partial_t u(x,t) &=& (H_{\Lambda} u)(x,t), &x \in \Lambda,\,t>0,\\
u(x,0) &=& \delta_y(x), &x \in \Lambda,
\end{array}
\end{equation}
and set $U^y_\Lambda(t) := \sum_{x \in \Lambda} u^y_\Lambda(x,t)$. Let $\tau_{\Lambda^\cc}$ be the hitting time of $\Lambda^\cc$ and 
\begin{equation}
\label{e:FKformula}
u^y_\Lambda(x,t) = \E_y \left[ \exp \left\{\int_0^t \xi(X_s) \textd s \right\} 
\1 \{\tau_{\Lambda^{\cc}}>t, X_t = x\} \right],
\end{equation}
the Feynman-Kac solution to \eqref{e:PAMdef} with Dirichlet boundary conditions on $\Lambda^\cc$. Then $u^y_\Lambda(x,t) $ also admits the spectral representation
\begin{equation}
\label{e:specrepr}
u^y_\Lambda(x,t) = \sum_{k=1}^{|\Lambda|} \texte^{t \lambda^{\ssup k}_\Lambda} 
\phi_\Lambda^{\ssup k}(y) \phi^{\ssup k}_\Lambda(x),
\end{equation}
where $\lambda^{\ssup 1}_\Lambda \ge \lambda^{\ssup 2}_\Lambda \ge \cdots \ge \lambda^{\ssup{|\Lambda|}}_\Lambda$ and $\phi^{\ssup 1}_\Lambda, \phi^{\ssup 2}_\Lambda, \ldots, \phi^{\ssup{|\Lambda|}}_\Lambda$ are, respectively the eigenvalues and the corresponding orthonormal eigenfunctions of $\Delta +\xi$ restricted to $\Lambda$. These two representations may be exploited to obtain bounds for one in terms of the other, as shown by the following lemma.
\begin{lemma}{\bf [Bounds on the solution]}
\label{l:bounds_mass}
For any $y \in \Lambda$ and any $t > 0$,
\begin{equation}
\label{e:bounds_mass}
\texte^{t \lambda^{\ssup 1}_\Lambda} \phi^{\ssup 1}_\Lambda(y)^2 
\le \E_y \left[ \texte^{\int_0^t \xi(X_s) \textd s} \1_{\{\tau_{\Lambda^\cc} > t, X_t = y\}} \right] \\
\le \E_y \left[ \texte^{\int_0^t \xi(X_s) \textd s} \1_{\{\tau_{\Lambda^\cc} > t\}} \right]
\end{equation}
\end{lemma}

\begin{proof}
The first inequality follows from a suitable application of Parseval's identity. The second inequality is elementary.
\end{proof}
\begin{lemma}{\bf [Mass up to an exit time]}
\label{l:mass_out}
For any $y \in \Lambda$, $\xi \in [0,\infty)^V$ and $\gamma > \lambda_\Lambda = \lambda_\Lambda(\xi,\GW)$,
\begin{equation}
\label{e:mass_out}
\EE_y \left[\ee^{\int_0^{\tau_{\Lambda^\cc}} (\xi(X_s) - \gamma )\, \textd s} \right] 
\le 1 + \frac{|\Lambda|}{ \gamma  - \lambda_\Lambda}.
\end{equation}
\end{lemma}
\begin{proof}
We follow the proof of \cite[Lemma 3.2]{dhW2022} and \cite[Lemma 2.18]{GM1998}. Define 
\begin{equation}
\label{e:u}
u(x) := \EE_x \left[\ee^{ \int_0^{\tau_{\Lambda^\cc}} (\xi(X_s) - \gamma )\, \textd s} \right].
\end{equation}
This is the solution to the boundary value problem
\begin{align}
\begin{split}
(\Delta + \xi -\gamma)u &= 0 \quad \text{on}\ \Lambda,\\
u &=1 \quad \text{on}\ \Lambda^\cc.
\end{split}
\end{align}
Via the substitution $u=:1+v$, this turns into 
\begin{align}
\begin{split}
(\Delta + \xi -\gamma)v &= \gamma - \xi \quad \text{on}\ \Lambda,\\
v &=0 \qquad\quad \text{on}\ \Lambda^\cc. 
\end{split}
\end{align}
It is readily checked that for $\gamma  > \lambda_\Lambda$ the solution exists and is given by
\begin{equation}
v = \mathcal{R}_\gamma(\xi-\gamma),
\end{equation}
where $\mathcal{R}_\gamma$ denotes the resolvent of $\Delta + \xi$. Hence
\begin{equation}
v(x) \leq (\mathcal{R}_\gamma\mathds{1})(x)
\leq \langle \mathcal{R}_\gamma\mathds{1},\mathds{1}\rangle
\leq \frac{|\Lambda|}{\gamma-\lambda_\Lambda}, \quad x \in \Lambda,  
\end{equation}
where $\mathds{1}$ denotes the constant function equal to $1$, and $\langle\cdot,\cdot\rangle$ denotes the weighted inner product. To get the first inequality, we apply the lower bound in \eqref{e:monot_princev} from Lemma~\ref{lem:specbd}, to get $\xi - \gamma \leq \lambda_\Lambda + 1 -\gamma \leq 1$ on $\Lambda$. The positivity of the resolvent gives
\begin{equation}
0 \leq [\mathcal{R}_\gamma(1 - (\xi-\gamma))](x) 
= [\mathcal{R}_\gamma\mathds{1}](x) - [\mathcal{R}_\gamma(\xi-\gamma)](x).
\end{equation}
To get the second inequality, we write
\begin{equation}
(\mathcal{R}_\gamma\mathds{1})(x) \leq \sum_{x \in \Lambda} (\mathcal{R}_\gamma\mathds{1})(x) 
=  \sum_{x \in \Lambda} (\mathcal{R}_\gamma\mathds{1})(x)\mathds{1}(x) \leq  \sum_{x \in \Lambda} (\mathcal{R}_\gamma\mathds{1})(x)\mathds{1}(x)\deg(x)
= \langle \mathcal{R}_\gamma\mathds{1},\mathds{1}\rangle.
\end{equation}
To get the third inequality, we use the Fourier expansion of the resolvent with respect to the orthonormal basis of eigenfunctions of $\Delta + \xi$ in $\langle \cdot, \cdot \rangle$.
\end{proof}
\subsection{Structural properties of the Galton-Watson tree}
\label{s:struct}
All of the results below can be lifted directly from \cite{dhW2022} since the normalisation of the Laplacian plays no role for the properties of the Galton-Watson tree. The results are included for the sake of completeness. 
\noindent
\begin{lemma}{\bf [Maximal degree in a ball around the root]}
\label{lem:deginball}
$\mbox{}$\\
(a) Subject to Assumption~{\rm \ref{ass:deg}(2)}, for every $\delta>0$,
\begin{equation}
\sum_{r \in \mathbb{N}} \mathfrak{P}\big(\exists\,x \in B_{2r}(\cO)\colon\, \deg(x) \geq  \delta r \big) < \infty.
\end{equation}
\end{lemma}
\begin{proof}
See \cite[Lemma 2.3]{dhW2022}.   
\end{proof}
\noindent Lemma~\ref{lem:deginball} shows that $\mathfrak{P}$-almost surely, as $r \to \infty$ all degrees in a ball of radius $r$ are eventually less than $\delta r$ for any $\delta >0$.
\begin{lemma}{\bf [Volumes of large balls]}
\label{lem:vol}
If there exists an $a>0$ such that $\mathcal{E}[\ee^{aD}] <\infty$, then for any $R_r$ satisfying $\lim_{r\to\infty} R_r/\log r = \infty$, 
\begin{equation}
\label{growth}
\liminf_{r\to\infty} \frac{1}{R_r} \log \Big(\inf_{x \in B_r(\cO)} |B_{R_r}(x)|\Big) 
= \limsup_{r\to\infty} \frac{1}{R_r} \log \Big(\sup_{x \in B_r(\cO)} |B_{R_r}(x)|\Big) 
= \vartheta \qquad \mathfrak{P}-a.s.
\end{equation} 
\end{lemma}
\begin{proof}
See \cite[Lemma 2.2]{dhW2022}.
\end{proof}
\noindent Lemmma~\ref{lem:vol} gives that $\mathfrak{P}$-almost surely, any ball 
of radius $r$ centred within distance $r$ to the root also has volume $\ee^{r\vartheta + o(1)}$ as $r \to \infty$. 
\subsection{Estimates on the potential}
\label{ss:potestim}
All of the results below are lifted directly from \cite{dhW2022}, since the normalisation of the Laplacian plays no role in the properties of the potential. The results are included for the sake of completeness. 
Abbreviate $L_r = |B_r(\cO)|$ and put
\begin{equation}
\label{e:def_Sr}
S_r := (\log r)^\alpha, \qquad \alpha \in (0,1).
\end{equation} 
For every $r \in \N$ there is a unique $a_r$ such that
\begin{equation}
\Prob(\xi(0) > a_r) = \frac{1}{r}.
\end{equation}
By Assumption~\ref{ass:pot}, for $r$ large enough
\begin{equation}
\label{def:Alr}
a_r = \varrho \log\log r.
\end{equation}
For $r \in \N$ and $A>0$, let
\begin{equation}
\label{def:Pi}
\Pi_{r,A} = \Pi_{r,A}(\xi) := \{z \in B_r(\cO)\colon\,\xi(z)>a_{L_r}-2A\}
\end{equation}
be the set of vertices in $B_r(\cO)$ where the potential is close to maximal,
\begin{equation}
\label{def:D}
D_{r,A} = D_{r,A}(\xi) := \{z \in B_r(\cO)\colon\,\mathrm{dist}(z,\Pi_{r,A}) \leq S_r\}
\end{equation}
be the $S_r$-neighbourhood of $\Pi_{r,A}$, and $\mathfrak{C}_{r,A}$ be the set of connected components of $D_{r,A}$ in $\GW$, which we think of as \emph{islands}. For $M_A\in\N$, define the event
\begin{equation}
\cB_{r,A} := \big\{ \exists\, \cC \in \mathfrak{C}_{r,A}\colon\, |\cC \cap \Pi_{r,A}| > M_A \big\}.
\end{equation}
Note that $\Pi_{r,A}, D_{r,A}, \cB_{r,A}$ depend on $\GW$ and therefore are random.
\begin{lemma}{\bf [Maximum size of the islands]}
\label{lem:size}
Subject to Assumptions~{\rm \ref{ass:pot}--\ref{ass:deg}}, for every $A > 0$ there exists an $M_A \in \N$ such that 
\begin{equation}
\sum_{r \in \N} \Prob(\cB_{r,A}) < \infty \qquad \mathfrak{P}-a.s.
\end{equation}
\end{lemma}
\begin{proof}
See \cite[Lemma 3.1]{dhW2022} and \cite[Lemma 6.6]{BK2016}.
\end{proof}
\noindent Lemma~\ref{lem:size} implies that $(\Prob\times\mathfrak{P})$-a.s.\ $\cB_{r,A}$ does not occur eventually as $r \to \infty$. Note that $\mathfrak{P}$-a.s.\ on the event $[\cB_{r,A}]^c$, 
\begin{equation}
\label{sizeresults}
\forall\,\cC \in \mathfrak{C}_{r,A}\colon\,
|\cC \cap \Pi_{r,A}| \leq M_A, \, \diam_\GW(\cC) \leq 2M_A S_r, \, |\cC| \leq \eee^{2\vartheta M_AS_r},
\end{equation}
where the last inequality follows from Lemma~\ref{lem:vol}. 
\begin{lemma}{\bf [Maximum of the potential]}
\label{l:maxpotential}
Subject to Assumptions~{\rm \ref{ass:pot}--\ref{ass:deg}}, for any $\vartheta>0$, $(\Prob \times \mathfrak{P})$-a.s.\ eventually as $r \to \infty$,
\begin{equation}
\label{e:asmaxpot}
\left| \max_{x \in B_r(\cO)} \xi(x) - a_{L_r} \right| \leq \frac{2 \varrho \log r}{\vartheta r}.
\end{equation}
\end{lemma}
\begin{proof}
See \cite[Lemma 2.4]{dHKdS2020}. The proof carries over verbatim and uses Lemma~\ref{lem:vol}.
\end{proof}
\begin{lemma}{\bf [Number of intermediate peaks of the potential]}
\label{l:bound_mediumpoints}
Subject to Assumptions~{\rm \ref{ass:pot}} and {\rm \ref{ass:deg}(2)}, for any $\beta \in (0,1)$ and $\varepsilon \in (0, \tfrac12\beta)$ the following holds. For a self-avoiding path $\pi$ in $\GW$, set
\begin{equation}
\label{e:bound_mediumpoints}
N_{\pi} = N_{\pi}(\xi) :=|\{z \in \supp(\pi) \colon\, \xi(z) > (1-\varepsilon) a_{L_r} \}|.
\end{equation}
Define the event
\begin{equation}
\cB_r := \left\{ 
\substack{\text{there exists a self-avoiding path } \pi \text{ in $\GW$ with } \\ 
 \supp(\pi) \cap B_r(\cO)\neq \emptyset, \, |\supp(\pi)| \geq (\log L_r)^{\beta}
\text{ and }  N_\pi > \frac{|\supp(\pi)|}{(\log {L_r})^\varepsilon}}
\right\}.
\end{equation}
Then  
\begin{equation}
\sum_{r \in \N_0} \Prob(\cB_r) < \infty \qquad \mathfrak{P}-a.s.
\end{equation}
\end{lemma}
\begin{proof}
See \cite[Lemma 3.6]{dhW2022}.    
\end{proof}
\noindent Lemma~\ref{l:bound_mediumpoints} implies that $(\Prob \times \mathfrak{P})$-a.s.\ for $r$ large enough, all self-avoiding paths $\pi$ in $\GW$ with $\supp(\pi) \cap B_r(\cO)\neq \emptyset$ and $|\supp(\pi)| \geq (\log L_r)^{\beta}$ satisfy $N_{\pi} \le \frac{|\supp(\pi)|}{(\log L_r)^\varepsilon}$.
\begin{lemma}{\bf [Number of high exceedances of the potential]}
\label{l:boundhighexceedances}
Subject to Assumptions~{\rm \ref{ass:pot}} and {\rm \ref{ass:deg}(2)}, for any $A>0$ there is a $C \ge 1$ such that, for all $\delta \in (0,1)$, the following holds. For a self-avoiding path $\pi$ in $\GW$, let
\begin{equation}
N_\pi := |\{ x \in \supp(\pi) \colon\, \xi(x) > a_{L_r} - 2A \}|.
\end{equation}
Define the event
\begin{equation}
\BB_r := \left\{ 
\substack{\text{there exists a self-avoiding path } \pi \text{ in $G$ with } \\ 
 \supp(\pi) \cap B_r(\cO)\neq \emptyset, \, |\supp(\pi)| \geq C (\log L_r)^{\delta}
\text{ and }  N_\pi > \frac{|\supp(\pi)|}{(\log {L_r})^\delta}}
\right\}.
\end{equation}
Then $\sum_{r \in \N_0} \sup_{G \in \mathfrak{G}_r} \Prob(\BB_r) < \infty$. In particular, $(\Prob \times \mathfrak{P})$-a.s.\ for $r$ large enough, all self-avoiding paths $\pi$ in $\GW$ with $\supp(\pi) \cap B_r(\cO)\neq \emptyset$ and $|\supp(\pi)| \geq C (\log L_r)^{\delta}$ satisfy
\begin{equation}
\label{e:boundhighexceedances}
N_\pi = |\{ x \in \supp(\pi) \colon\, \xi(x) > a_{L_r} - 2A \}| \le \frac{|\supp(\pi)|}{(\log L_r)^\delta}.
\end{equation} 
\end{lemma}
\begin{proof}
See \cite[Lemma 3.7]{dhW2022}.
\end{proof}
\begin{lemma}{\bf [Principal eigenvalues of the islands]}
\label{l:eigislands}
Subject to Assumptions~{\rm \ref{ass:pot}} and {\rm \ref{ass:deg}(2)}, for any $\varepsilon>0$, $(\Prob\times\mathfrak{P})$-a.s.\ eventually as $r \to \infty$,
\begin{equation}
\label{e:eigislands}
\text{all}\ \cC \in \mathfrak{C}_{r,A}\ \text{satisfy} \colon\, \lambda_\cC(\xi; \GW) 
\leq a_{L_r} - \widehat{\chi}_{\cC}(\GW) + \varepsilon.
\end{equation}
\end{lemma}
\begin{proof}
See \cite[Lemma 3.3]{dhW2022}.
\end{proof}
\begin{corollary}{\bf [Uniform bound on principal eigenvalue of the islands]}
\label{c:eigislandsGW}
Subject to Assumptions~{\rm \ref{ass:pot}--\ref{ass:deg}}, for $\vartheta$ as in \eqref{e:volumerateGW}, and any $\varepsilon>0$, $(\Prob \times \Probgr)$-a.s.\ eventually as $r \to \infty$,
\begin{equation}
\label{e:eigislandsGW}
\max_{\CC \in \mathfrak{C}_{r,A}}\lambda^{\ssup 1}_\CC(\xi; G) 
\leq a_{L_r} - \widetilde{\chi}(\varrho) + \varepsilon.
\end{equation}
\end{corollary}
\begin{proof}
See \cite[Corollary 2.8]{dHKdS2020}.    
\end{proof}

\section{Path expansions}
\label{s.pathexpansions}
In this section we adapt \cite[Section 3]{dHKdS2020} to fit with the random walk generated by the normalised Laplacian. Section~\ref{preplem} proves three lemmas that concern the contribution to the total mass in \eqref{e:mass} coming from various sets of paths. Section~\ref{ss:keyprop} proves a key proposition that controls the entropy associated with a key set of paths. The proof of which is based on the three lemmas in Section~\ref{preplem}.

We need various sets of nearest-neighbour paths in $\GW=(V,E,\cO)$, defined in \cite{dHKdS2020}. For $\ell \in \N_0$ and subsets $\Lambda, \Lambda' \subset V$, put
\begin{equation}
\label{defcurlyP}
\begin{aligned}
&\scrP_\ell(\Lambda,\Lambda') := \left\{ (\pi_0, \ldots, \pi_{\ell}) \in V^{\ell+1} \colon\,
\begin{array}{ll} 
&\pi_0 \in \Lambda, \pi_{\ell} \in \Lambda',\\
&\{\pi_{i}, \pi_{i-1}\} \in E \;\forall\, 1 \le i \le \ell
\end{array}
\right\},\\
&\scrP(\Lambda, \Lambda') := \bigcup_{\ell \in \N_0} \scrP_\ell(\Lambda,\Lambda'),
\end{aligned}
\end{equation}
and set
\begin{equation}
\scrP_\ell := \scrP_\ell(V,V), \qquad \scrP := \scrP(V,V). 
\end{equation}
When $\Lambda$ or $\Lambda'$ consists of a single point, write $x$ instead of $\{x\}$. For $\pi \in \scrP_\ell$, set $|\pi| := \ell$. Write $\supp(\pi) := \{\pi_0, \ldots, \pi_{|\pi|}\}$ to denote the set of points visited by $\pi$.

Let $X=(X_t)_{t\ge0}$ be the continuous-time random walk on $G$ that jumps from $x \in V$ to any neighbour $y\sim x$ at rate $1$. Denote by $(T_k)_{k \in \N_0}$ the sequence of jump times (with $T_0 := 0$). For $\ell \in \N_0$, let 
\begin{equation}
\pi^{\ssup \ell}(X) := (X_0, \ldots, X_{T_{\ell}})
\end{equation}
be the path in $\scrP_\ell$ consisting of the first $\ell$ steps of $X$. For $t  \ge 0$, let
\begin{equation}
\label{e:defpathX0t}
\pi(X_{[0,t]}) = \pi^{\ssup{\ell_t}}(X), \quad \text{ with } \ell_t \in \N_0 \, 
\text{ satisfying } \, T_{\ell_t} \le t < T_{\ell_t+1},
\end{equation}
denote the path in $\scrP$ consisting of all the steps taken by $X$ between times $0$ and $t$.

Recall the definitions from Section~\ref{ss:potestim}. For $\pi \in \scrP$ and $A>0$, define
\begin{equation}
\label{e:deflambdaLApi}
\lambda_{r,A}(\pi) := \sup \big\{ \lambda^{\ssup 1}_\CC(\xi; G) 
\colon\, \CC \in \mathfrak{C}_{r,A}, 
\, \supp(\pi)\cap \CC \cap \Pi_{r,A} \neq \emptyset \big\},
\end{equation}
with the convention $\sup \emptyset = -\infty$. 
This is the largest principal eigenvalue among the components of $\mathfrak C_{r,A}$ in $\GW$ 
that have a point of high exceedance visited by the path $\pi$.

\subsection{Mass of the solution along excursions}
\label{preplem}

\begin{lemma}{\bf [Path evaluation]}
\label{l:path_eval}
For $\ell\in\N_0$, $\pi \in \scrP_\ell$ and $\gamma  > \max_{0 \leq i < |\pi|} \{\xi(\pi_i)-1\}$,
\begin{equation}
\label{e:path_eval}
\E_{\pi_0} \left[\ee^{\int_0^{T_{\ell}} (\xi(X_s) -  \gamma )\, \textd s} ~\Big|~ \pi^{\ssup {\ell}}(X) = \pi  \right]
= \prod_{i=0}^{\ell-1} \frac{1}{\gamma - [\xi(\pi_i)-1]}.
\end{equation}
\end{lemma}

\begin{proof}
The proof is identical to that of \cite[Lemma 3.2]{dHKdS2020}, except that the random walk now jumps with rate 1.
\end{proof}

For a path $\pi \in \scrP$ and $\varepsilon \in (0,1)$, we write
\begin{equation}
\label{e:def_Mpi}
M^{r,\varepsilon}_\pi := \big| \bigl\{0 \leq i < |\pi| \colon\, \xi(\pi_i) \le (1-\varepsilon)a_{L_r}\bigr\}\big|,
\end{equation}
with the interpretation that $M^{r,\varepsilon}_\pi = 0$ if $|\pi|=0$.

\begin{lemma}{\bf [Mass of excursions]}
\label{l:mass_in}
Subject to Assumption~{\rm \ref{ass:pot}}, for every $A, \varepsilon>0$, there exists $c>0$ and $r_0 \in \N$ such that, for all $r \ge r_0$, all $\gamma > a_{L_r} - A$ and all $\pi \in \scrP(B_r(\cO), B_r(\cO))$ satisfying $\pi_i \notin \Pi_{r,A}$ for all $0 \leq i < \ell:=|\pi|$,
\begin{equation}
\label{e:mass_in}
\E_{\pi_0} \left[\ee^{ \int_0^{T_{\ell}}(\xi(X_s) - \gamma)\, \textd s} ~\Big|~ \pi^{\ssup {\ell}}(X) = \pi \right]
\leq q_{r,A}^{\ell} \ee^{ M^{r,\varepsilon}_\pi (c-\log\log\log L_r)},
\end{equation}
where 
\begin{equation}
\label{Aqdef} \qquad q_A := \frac{1}{1+A} \quad \mathrm{and} \quad c=\log[2(q_A\varepsilon\varrho)^{-1}].
\end{equation} 
Note that $\pi_{\ell} \in \Pi_{r,A}$ is allowed.
\end{lemma}

\begin{proof}
The proof is identical to that of \cite[Lemma 3.3]{dHKdS2020}, and uses Lemma~\ref{l:path_eval}. 
\end{proof}

We follow \cite[Definition 3.4]{dHKdS2020} and \cite[Section 6.2]{BKS2018}. Note that the distance between $\Pi_{r,A}$ and $D_{r,A}^\cc$ in $\GW$ is at least $S_r = (\log L_r)^\alpha$ (recall \eqref{def:Pi}--\eqref{def:D}).

\begin{definition}{\bf [Concatenation of paths]} {\rm (a)}
\label{def:concat}
When $\pi$ and $\pi'$ are two paths in $\scrP$ with $\pi_{|\pi|} = \pi'_0$, 
we define their \emph{concatenation} as
\begin{equation}
\label{def_concat}
\pi \circ \pi' := (\pi_0, \ldots, \pi_{|\pi|}, \pi'_1, \ldots, \pi'_{|\pi'|}) \in \scrP.
\end{equation}
Note that $|\pi \circ \pi'| = |\pi| + |\pi'|$. 

\medskip\noindent
{\rm (b)} When $\pi_{|\pi|} \neq \pi'_0$, we can still define the \emph{shifted concatenation} of $\pi$ and $\pi'$ as $\pi \circ \hat{\pi}'$, where $\hat{\pi}' := (\pi_{|\pi|}, \pi_{|\pi|}  + \pi'_1 - \pi'_0, \ldots, \pi_{|\pi|} + \pi'_{|\pi'|} - \pi'_0)$. The shifted concatenation of multiple paths is defined inductively via associativity. 
\end{definition}

Now, if a path $\pi \in \scrP$ intersects $\Pi_{r,A}$, then it can be decomposed into an initial path, 
a sequence of excursions between $\Pi_{r,A}$ and $D_{r,A}^\cc$, and a terminal path. 
More precisely, there exists $m_\pi \in \N $ such that
\begin{equation}
\label{e:concat1}
\pi = \check{\pi}^{\ssup 1} \circ \hat{\pi}^{\ssup 1} \circ \cdots \circ \check{\pi}^{\ssup {m_\pi}} 
\circ \hat{\pi}^{\ssup {m_\pi}} \circ \bar{\pi},
\end{equation}
where the paths in \eqref{e:concat1} satisfy
\begin{equation}
\label{e:concat2}
\begin{alignedat}{9}
\check{\pi}^{\ssup 1} & \in  \scrP(V, \Pi_{r,A}) 
&\qquad\text{with}\qquad& 
\check{\pi}^{\ssup 1}_i & \notin  \Pi_{r,A}, & \quad\, 0\le i < |\check{\pi}^{\ssup 1}|, 
\\
\hat{\pi}^{\ssup k} & \in  \scrP(\Pi_{r,A}, D_{r,A}^\cc) 
&\qquad\text{with}\qquad& 
\hat{\pi}^{\ssup k}_i & \in  D_{r, A}, & \quad\, 0\le i < |\hat{\pi}^{\ssup k}|, \; 1 \le k \le m_{\pi} - 1, 
\\
\check{\pi}^{\ssup k} & \in  \scrP(D_{r,A}^\cc, \Pi_{r,A}) 
&\qquad\text{with}\qquad& 
\check{\pi}^{\ssup k}_i & \notin  \Pi_{r,A}, & \quad\, 0\le i < |\check{\pi}^{\ssup k}|, \; 2 \le k \le m_\pi, 
\\
\hat{\pi}^{\ssup {m_\pi}} & \in  \scrP(\Pi_{r,A}, V) 
&\qquad\text{with}\qquad& 
\hat{\pi}^{\ssup {m_\pi}}_i & \in  D_{r,A}, & \quad\, 0\le i < |\hat{\pi}^{\ssup {m_\pi}}|, 
\end{alignedat}
\end{equation}
while
\begin{equation}
\label{e:concat3}
\begin{array}{ll} 
\bar{\pi} \in \scrP(D_{r,A}^\cc, V) \text{ and } \bar{\pi}_i \notin \Pi_{r,A} \; \forall\, i \ge 0 
& \text{ if } \hat{\pi}^{\ssup {m_\pi}} \in \scrP(\Pi_{r,A}, D^\cc_{r, A}), \\
\bar{\pi}_0 \in D_{r,A}, |\bar{\pi}| = 0  & \text{ otherwise.}
\end{array}
\end{equation}
Note that the decomposition in \eqref{e:concat1}--\eqref{e:concat3} is unique, and that the paths $\check{\pi}^{\ssup 1}$, $\hat{\pi}^{\ssup {m_\pi}}$ and $\bar{\pi}$ can have zero length. If $\pi$ is contained in $B_r(\cO)$, then so are all the paths in the decomposition. 

Whenever $\supp(\pi) \cap \Pi_{r,A} \ne \emptyset$ and $\varepsilon > 0$, we define
\begin{align}
\label{e:defnpikpi}
s_\pi := \sum_{i=1}^{m_\pi} |\check{\pi}^{\ssup i}| + |\bar{\pi}|, \qquad
k^{r,\varepsilon}_\pi := \sum_{i=1}^{m_\pi} M^{r,\varepsilon}_{\check{\pi}^{\ssup i}} + M^{r,\varepsilon}_{\bar{\pi}},
\end{align}
to be the total time spent in exterior excursions, respectively, on moderately low points of the potential visited by exterior excursions (without their last point). 

In case $\supp(\pi) \cap \Pi_{r,A} = \emptyset$, we set $m_\pi := 0$, $s_\pi := |\pi|$ and $k^{r,\varepsilon}_\pi := M^{r,\varepsilon}_{\pi}$. Recall from \eqref{e:deflambdaLApi} that, in this case, $\lambda_{r,A}(\pi) = -\infty$. 

We say that $\pi, \pi' \in \scrP$ are \emph{equivalent}, written $\pi' \sim \pi$, if $m_{\pi} = m_{\pi'}$, $\check{\pi}'^{\ssup i}=\check{\pi}^{\ssup i}$ for all $i=1,\ldots,m_{\pi}$, and $\bar{\pi}' = \bar{\pi}$. If $\pi' \sim \pi$, then $s_{\pi'}$, $k^{r, \varepsilon}_{\pi'}$ and $\lambda_{r,A}(\pi')$ are all equal to the counterparts for $\pi$.

To state our key lemma, we define, for $m,s \in \N_0$,
\begin{equation}
\label{e:defPmn}
\scrP^{(m,s)} = \left\{ \pi \in \scrP \colon\, m_\pi = m, s_\pi = s \right\},
\end{equation}
and denote by
\begin{equation}
\label{def_CLA}
C_{r,A}:= \max \{|\CC| \colon\, \CC \in \mathfrak{C}_{r,A}\}
\end{equation}
the maximal size of the islands in $\mathfrak{C}_{r,A}$.

\begin{lemma}{\bf [Mass of an equivalence class]}
\label{l:fixed_class}
For every $A,\varepsilon > 0$ there exist $c>0$ and $r_0 \in \N$ such that, for all $r \ge r_0$, all $m,s \in \N_0$, all $\pi \in \scrP^{(m,s)}$ with $\supp(\pi) \subset B_r(\cO)$, all $\gamma > \lambda_{r,A}(\pi) \vee (a_{L_r} -A)$ and all $t \ge 0$,    
\begin{multline}
\label{e:fixed_class}
\qquad
\E_{\pi_0} \left[ \texte^{\int_0^t (\xi(X_u) - \gamma)\, \textd u}\, \1_{\{\pi(X_{[0,t]}) \sim \pi\}} \right]  
\\\le \left(C_{r,A}^{1/2} \right)^{\1_{\{m>0\}}} \left(1+\frac{ C_{r,A}}{\gamma - \lambda_{r,A}(\pi)} \right)^m \left(\frac{q_A}{d_{\min}}\right)^s \texte^{\left(c-\log^{\ssup 3} L_r \right) k^{r,\varepsilon}_{\pi}}.
\end{multline}
\end{lemma}
\begin{proof}
The proof is identical to that of \cite[Lemma 3.5]{dHKdS2020}, except that the normalised Laplacian gives rise to Lemma~\ref{l:mass_out} and Lemma~\ref{l:mass_in}, which are used instead.
\end{proof}


\subsection{Key proposition}
\label{ss:keyprop}

The main result of this section is the following proposition.

\begin{proposition}{\bf [Entropy reduction]}
\label{p:massclass}
Let $\alpha \in (0,1)$ be as in \eqref{e:def_Sr} and $\kappa\in (\alpha,1)$. Subject to Assumption~{\rm \ref{ass:deg}}, there exists an $A_0$ such that, for all $A \geq A_0$, with $\mathfrak{P}$-probability tending to one as $r\to\infty$, the following statement is true. For each $x \in B_r(\cO)$, each $\NN \subset \scrP(x,B_r(\cO))$ satisfying $\supp(\pi) \subset B_r(\cO)$ and $\max_{1 \le \ell \le |\pi|} \dist_{G}(\pi_\ell, x) \geq (\log L_r)^\kappa$ for all $\pi \in \mathcal{N}$, and each assignment $\pi\mapsto (\gamma_\pi , z_\pi)\in \R \times V$ satisfying
\begin{equation}
\label{e:cond_massclass1}
\gamma_\pi \ge \left(\lambda_{r,A}(\pi)  + \texte^{-S_r} \right) \vee (a_{L_r}- A) \qquad \forall\,\,\pi \in \NN
\end{equation}
and
\begin{equation}
\label{e:cond_massclass2}
z_\pi \in \supp(\pi) \cup 
\bigcup_{ \substack{\CC \in \mathfrak{C}_{r,A} \colon \\ \supp(\pi) \cap \CC \cap \Pi_{r,A} \neq \emptyset}} \CC 
\qquad \forall\,\, \pi \in \NN,
\end{equation}
the following inequality holds for all $t \geq 0$:
\begin{equation}
\label{e:mass_class}
\log \E_x \left[ \texte^{\int_0^t \xi(X_s) \textd s} \1_{\{\pi(X_{[0,t]}) \in \mathcal{N}\}}\right]
\leq \sup_{\pi \in \mathcal{N}} \Big\{ t \gamma_\pi + \dist_{G}(x,z_\pi) (c-\log\log\log L_r)\Big\}. 
\end{equation}
\end{proposition}

\begin{proof}
The proof is based on \cite[Section 3.4]{dHKdS2020}. First fix $c_0 >2$  and define
\begin{equation}
\label{e:defc0A0}
A_0 =  \ee^{4 c_0}-1.
\end{equation}
Fix $A \geq A_0$, $\beta \in (0,\alpha)$ and $\varepsilon \in (0,\frac12\beta)$ as in Lemma~\ref{l:bound_mediumpoints}. Let $r_0 \in \N$ be as given in Lemma~\ref{l:fixed_class}, and take $r \ge r_0$ so large that the conclusions of Lemmas~\ref{lem:deginball}, \ref{lem:size}, \ref{l:eigislands} and \ref{l:bound_mediumpoints} hold, i.e. assume that the events $\BB_r$ and $\BB_{r,A}$ in these lemmas do not occur. Fix $x \in B_r(\cO)$. Recall the definitions of $C_{r,A}$ and $\scrP^{(m,s)}$. Note that the relation $\sim$ is an equivalence relation in $\scrP^{(m,s)}$, and define
\begin{equation}
\label{e:propmassclass2}
\widetilde{\scrP}^{(m,s)}_x := \big\{\text{equivalence classes of the paths in } \scrP(x,V) \cap \scrP^{(m,s)}\big\}.
\end{equation}
The following bounded on the cardinality of this set is needed.

\begin{lemma}{\bf [Bound equivalence classes]}
\label{l:propmassclass3}
Subject to Assumption~{\rm \ref{ass:deg}}, $\mathfrak{P}$-a.s.,$|\widetilde{\scrP}^{(m,s)}_x| $ $\le (2C_{r,A})^m (\delta r)^{(m+s)}$ for all $m,s \in \N_0$.
\end{lemma}

\begin{proof}
We can copy the proof of \cite[Lemma 3.6]{dHKdS2020}, replacing $d_{\max}$ by $\delta r$.
\end{proof}

Now take $\NN \subset \scrP(x, V)$ as in the statement, and set
\begin{equation}
\label{e:propmassclass1}
\widetilde{\mathcal{N}}^{(m,s)} := \big\{\text{equivalence classes of paths in } 
\NN \cap \scrP^{(m,s)}\big\} \subset \widetilde{\scrP}^{(m,s)}_x.
\end{equation}
For each $\MM \in \widetilde{\NN}^{(m,s)}$, choose a representative $\pi_\MM \in \MM$, and use Lemma~\ref{l:propmassclass3} to write 
\begin{align}
\label{e:propmassclass6}
& \E_x \left[ \texte^{\int_0^t \xi(X_u) \textd u} \1_{\{\pi(X_{[0,t]}) \in \mathcal{N}\}} \right] 
= \sum_{m, s \in \N_0}  \sum_{\MM \in \widetilde{\mathcal{N}}^{(m,s)}}
\E_x \left[ \texte^{\int_0^t \xi(X_u) \textd u} \1_{\{\pi(X_{[0,t]}) \sim \pi_\MM \}} \right] \nonumber\\
& \quad\qquad \le \sum_{m, s \in \N_0} (2 (\delta r) C_{r,A})^m (\delta r)^s 
\sup_{\pi \in \NN^{(m,s)}} \E_x \left[ \texte^{\int_0^t \xi(X_u) \textd u} \1_{\{\pi(X_{[0,t]}) \sim \pi\}} \right]
\end{align}
with the convention $\sup \emptyset = 0$. For fixed $\pi \in \mathcal{N}^{(m,s)}$, by \eqref{e:cond_massclass1}, apply \eqref{e:fixed_class} and Lemma~\ref{lem:size} to obtain, for all $r$ large enough and with $c_0 >2$ ,
\begin{equation}
\label{e:propmassclass7}
\begin{aligned}
&(2(\delta r))^m  (\delta r)^s\, 
\E_x \left[ \texte^{\int_0^t \xi(X_u) \textd u} \1_{\{\pi(X_{[0,t]}) \sim \pi\}} \right]\\ 
&\qquad \le \texte^{t \gamma_\pi } \texte^{c_0 m S_r} [q_A (\delta r)]^s\, 
\ee^{k^{r,\varepsilon}_\pi (c-\log\log\log L_r)}.
\end{aligned}
\end{equation}
We next claim that, for $r$ large enough and $\pi \in \NN^{(m,s)}$,
\begin{equation}
\label{e:propmassclass7.1}
s \ge \left[(m-1)\vee 1 \right] S_r .
\end{equation}
Indeed, when $m\ge 2$, $|\supp(\check{\pi}^{\ssup i})| \ge S_r$ for all $2 \le i \le m$. When $m=0$, $|\supp(\pi)| \ge \max_{1 \le \ell \le |\pi|} |\pi_\ell -x| \ge (\log L_r)^\kappa \gg S_r$ by assumption. When $m=1$, the latter assumption and Lemma~\ref{lem:size} together imply that $\supp(\pi) \cap D^\cc_{r,A} \neq \emptyset$, and so either $|\supp(\check{\pi}^{\ssup 1})| \geq S_r$ or $|\supp(\bar{\pi})|\ge S_r$. Thus, \eqref{e:propmassclass7.1} holds by the definition of $S_r$ and $s$.

Note that $q_A < \ee^{-4c_0}$, so
\begin{equation}
\label{e:propmassclass7.5}
\begin{aligned}
&\sum_{m \geq 0} \sum_{s \geq [(m-1)\vee 1] S_r} \ee^{c_0 m S_r} [q_A(\delta r)]^s
\\
&= \frac{[q_A (\delta r)]^{S_r} + \ee^{c_0 S_r}[q_A (\delta r)]^{S_r} + \sum_{m \geq 2} \ee^{mc_0 S_r  } [q_A (\delta r)]^{(m-1)S_r}}{1-q_A \delta r}\\
&\leq \frac{3 \ee^{-c_0 \log r}}{1-q_A \delta r} < 1
\end{aligned}
\end{equation}
for $r$ large enough. Inserting this back into~\eqref{e:propmassclass6}, we obtain
\begin{equation}
\label{e:intermediatemassclass}
\log \E_x \left[ \texte^{\int_0^t \xi(X_s) \textd s} \1_{\{\pi(X_{0,t}) \in \mathcal{N}\}} \right]
\leq \sup_{\pi \in \mathcal{N}} \Big\{ t \gamma_\pi 
+ k^{r,\varepsilon}_\pi (c-\log\log\log L_r)\Big\}.
\end{equation}
The remainder of the proof is identical to the end of \cite[Section 3.4]{dHKdS2020} and is included for completeness.

The proof will be finished once we show that, for some $\varepsilon' > 0$ and whp, respectively, a.s.\ eventually as $r \to \infty$, 
\begin{equation}
\label{e:propmassclass9}
k^{r,\varepsilon}_\pi \ge \dist_{G}(x,z_{\pi})(1-2(\log L_r)^{-\varepsilon'}) \qquad \forall\,\pi \in \NN.
\end{equation}
For each $\pi \in \NN$ define an auxiliary path $\pi_\star$ as follows. First note that by using our assumptions we can find points $z', z'' \in \supp(\pi)$ (not necessarily distinct) such that 
\begin{equation}
\label{e:propmassclass10}
\dist_{G}(x,z') \geq (\log L_r)^\kappa, \qquad \dist_{G}(z'', z_\pi) \leq 2 M_A S_r,
\end{equation}
where the latter holds by \eqref{sizeresults}. Write  $\{z_1, z_2 \} = \{z', z''\}$ with $z_1$, $z_2$ ordered according to their hitting times by $\pi$, i.e. $\inf\{ \ell \colon \pi_\ell = z_1 \} \leq \inf\{\ell \colon \pi_\ell = z_2\}$. Define $\pi_e$ as the concatenation of the loop erasure of $\pi$ between $x$ and $z_1$ and the loop erasure of $\pi$ between $z_1$ and $z_2$. Since $\pi_e$ is the concatenation of two self-avoiding paths, it visits each point at most twice. Finally, define $\pi_\star \sim \pi_e$ by replacing the excursions of $\pi_e$ from $\Pi_{r,A}$ to $D_{r,A}^\cc$ by direct paths between the corresponding endpoints, i.e. replace each $\hat{\pi}_e^{\ssup i}$ by $|\hat{\pi}_e^{\ssup i}|=\ell_i$, $(\hat{\pi}_e^{\ssup i})_0 = x_i \in \Pi_{r,A}$, and $(\hat{\pi}_e^{\ssup i})_{\ell_i} = y_i \in D_{r,A}^\cc$ by a shortest-distance path $\widetilde{\pi}_\star^{\ssup i}$ with the same endpoints and $|\widetilde{\pi}_\star^{\ssup i}| = \dist_{G}(x_i, y_i)$. Since $\pi_\star$ visits each $x \in \Pi_{r,A}$ at most $2$ times,
\begin{equation}
\label{e:propmassclass11}
\begin{aligned}
k^{r,\varepsilon}_\pi \ge k^{r,\varepsilon}_{\pi_\star} \geq M^{r,\varepsilon}_{\pi_\star} 
- 2 |\supp(\pi_\star)\cap \Pi_{r,A}|(S_r+1) \geq M^{r,\varepsilon}_{\pi_\star} - 4 |\supp(\pi_\star)\cap \Pi_{r,A}|  S_r.
\end{aligned}
\end{equation}
Note that $M_{\pi_\star}^{r, \varepsilon} \geq \left|\{x \in \supp(\pi_\star) \colon\, \xi(x) \leq (1-\varepsilon) a_{L_r}\} \right| - 1$ and, by \eqref{e:propmassclass10}, $|\supp(\pi_\star)| \geq \dist_{G}(x,z') \geq (\log L_r)^\kappa \gg (\log L_r)^{\alpha+2\varepsilon'}$ for some $0<\varepsilon'<\varepsilon$. Applying Lemmas~\ref{l:bound_mediumpoints}--\ref{l:boundhighexceedances} and using \eqref{e:def_Sr} and $L_r > r$, we obtain, for $r$ large enough,
\begin{equation}
\label{e:propmassclass12}
\begin{aligned}
k^{r,\varepsilon}_\pi 
& \geq |\supp(\pi_\star)|\left( 1 - \frac{2}{(\log L_r)^{\varepsilon}} 
-  \frac{4 S_r}{(\log L_r)^{\alpha+2\varepsilon'}}\right) 
\geq |\supp(\pi_\star)|\left( 1 - \frac{1}{(\log L_r)^{\varepsilon'}}\right).
\end{aligned}
\end{equation}
On the other hand, since $|\supp(\pi_\star)| \geq (\log L_r)^\kappa$, by \eqref{e:propmassclass10} we have
\begin{equation}
\label{e:propmassclass13}
\begin{aligned}
\left|\supp(\pi_\star) \right| &= \big(\left|\supp(\pi_\star) \right| +  2 M_A S_r\big) - 2 M_A S_r\\
&= \big(\left|\supp(\pi_\star) \right| +  2 M_A S_r\big) \left( 1- \frac{2 M_A S_r}{\left|\supp(\pi_\star) \right| +  2 M_A S_r}\right)\\
& \geq \left( \dist_{G}(x,z'') + 2 M_A S_r \right) \left( 1-\frac{2 M_A S_r}{(\log L_r)^\kappa} \right) \\
& \geq \dist_{G}(x,z_\pi)\left( 1-\frac{1}{(\log L_r)^{\varepsilon'}} \right),
\end{aligned}
\end{equation}
where the first inequality uses that the distance between two points on $\pi_\star$ is less than the total length of $\pi_\star$. Now \eqref{e:propmassclass9} follows from \eqref{e:propmassclass12}--\eqref{e:propmassclass13}.
\end{proof}


\section{Proof of the main theorem}
\label{s.proofmain}

Define
\begin{equation}
\label{e:ustar}
U^*(t) := \ee^{t[\varrho \log(\vartheta \mathfrak{r}_t) -\varrho - \widetilde{\chi}(\varrho)]},
\end{equation}
where we recall \eqref{mathfrakrdef}. To prove Theorem~\ref{t:QLyapGWT} we show that
\begin{equation}
\label{UU*comp}
\frac{1}{t} \log U(t) - \frac{1}{t} \log U^*(t) = o(1), \quad t \to \infty, 
\qquad (\Prob \times \Probgr)\text{-a.s.}
\end{equation}
The proof proceeds via upper and lower bound, proved in Sections~\ref{sec:ub} and \ref{sec:lb}, respectively.


\subsection{Upper bound}
\label{sec:ub}
We follow \cite[Section 4.2]{dHKdS2020}. The proof of the upper bound in \eqref{UU*comp} relies on two lemmas showing that paths staying inside a ball of radius $\lceil t^\gamma \rceil$ for some $\gamma \in (0,1)$ or leaving a ball of radius $t \log t$ have a negligible contribution to \eqref{e:mass}, the total mass of the solution.

\begin{lemma}{\bf [No long paths]}
\label{l:longpaths}
For any $\ell_t \geq t \log t$, 
\begin{equation}
\label{e:FT2}
\lim_{t \to \infty} \frac{1}{U^*(t)}\,\E_{\cO} \left[\ee^{\int_0^t \xi(X_s) \dd s} \1_{\{\tau_{[B_{\ell_t}]^\cc}< t\}}\right] = 0 
\quad (\Prob \times \mathfrak{P})-a.s.
\end{equation}
\end{lemma}
\begin{proof}
We follow \cite[Lemma 4.2]{dHKdS2020}. For $r \geq \ell_t$, let
\begin{equation}
\BB_r := \left\{ \max_{x \in B_r(\cO)} \xi(x) \geq a_{L_r} + 2 \varrho\right\}.
\end{equation}
Since $\lim_{t\to\infty} \ell_t = \infty$, Lemma~\ref{l:maxpotential} gives that $\Prob$-a.s.
\begin{equation}
\label{e:prFT2}
\bigcup_{r \geq \ell_t} \BB_r
\text{ does not occur eventually as } t\to \infty.
\end{equation}
Therefore we can work on the event $\bigcap_{r \geq \ell_t} [\BB_r]^\cc$. On this event, we write
\begin{align}
\label{e:prFT3}
\E_{\cO} \left[\ee^{\int_0^t \xi(X_s) \dd s} \1_{\{\tau_{[B_{\ell_t}]^\cc}< t\}} \right]
& = \sum_{r \geq \ell_t} \E_{\cO} \left[\ee^{\int_0^t \xi(X_s) \dd s} 
\1_{\{\sup_{s \in [0,t]}|X_s| = r \}} \right] \nonumber\\
& \leq  \ee^{2\varrho t} \sum_{r \geq \ell_t}\, \ee^{\varrho t \log r +  \varrho \log \log (\delta r)} \, 
\P_{\cO} \left( J_t \geq r \right),
\end{align}
where $J_t$ is the number of jumps of $X$ up to time $t$, and we use that $|B_r(\cO)| \leq (\delta r)^r$. Next, $J_t$ is stochastically dominated by a Poisson random variable with parameter $t$. Hence
\begin{equation}
\P_{\cO} \left( J_t \geq r \right) \leq \frac{(\ee t)^r}{r^r} \leq 
\exp \left\{-r \log\left( \frac{r}{\ee t} \right) \right\}
\end{equation}
for large $r$. Using that $\ell_t \geq  t \log t$, we can easily check that, for $r \geq \ell_t$ and $t$ large enough,
\begin{equation}
\varrho t \log r - r \log\left( \frac{r}{\ee t}\right) < -3 r, \qquad r \geq \ell_t. 
\end{equation}
Thus \eqref{e:prFT3} is at most
\begin{equation}
\ee^{2\varrho t} \sum_{r \geq \ell_t}\, \ee^{-3r + \log\log (\delta r)} \, \leq \ee^{2\varrho t} \sum_{r \geq \ell_t}\, \ee^{-2r} 
\leq 2\,\ee^{2\varrho t}\,\ee^{-2\ell_t} \leq \ee^{-\ell_t}. 
\end{equation}
Since $\lim_{t\to\infty} \ell_t = \infty$ and $\lim_{t\to\infty} U^*(t) = \infty$, this settles the claim.
\end{proof}

\begin{lemma}{\bf [No short paths]}
\label{l:noshortpaths}
For any $\gamma \in (0,1)$,
\begin{equation}
\label{e:NSP2}
\lim_{t \to \infty} \frac{1}{U^*(t)}\,\E_{\cO} \left[\ee^{\int_0^t \xi(X_s) \dd s} 
\1_{\{\tau_{[B_{\lceil t^\gamma \rceil}]^\cc} > t\}} \right] = 0 
\quad (\Prob \times \mathfrak{P})-a.s.
\end{equation}
\end{lemma}
\begin{proof}
We follow \cite[Lemma 4.3]{dHKdS2020}. By Lemma~\ref{l:maxpotential} with $r = \lceil t^\gamma \rceil$, we may assume that
\begin{equation}
\max_{x \in B_{\lceil t^\gamma \rceil}} \xi(x) \leq \varrho \log \log L_{\lceil t^\gamma \rceil}
+ \frac{2 \varrho \log\lceil t^\gamma \rceil}{\vartheta \lceil t^\gamma \rceil}
\leq \gamma \varrho \log t + O(1), \quad t \to \infty,
\end{equation}
where the second inequality uses that $\log L_{\lceil t^\gamma \rceil} \sim \log |B_{\lceil t^\gamma \rceil}(\cO)| \sim \vartheta \lceil t^\gamma \rceil$. Hence
\begin{equation}
\frac{1}{U^*(t)} \,\E_{\cO} \left[\ee^{\int_0^t \xi(X_s) \dd s} \1_{\{\tau_{[B_{\lceil t^\gamma \rceil}]^\cc} > t\}}\right]
\leq \frac{1}{U^*(t)}\,\ee^{\gamma\varrho t \log t+O(1)} \leq \ee^{ -(1-\gamma)\varrho t \log t + C t \log\log\log t}, \quad t \to \infty,
\end{equation}
for any constant $C>1$.
\end{proof}

The proof of the upper bound in \eqref{UU*comp} also relies on a third lemma estimating the contribution of paths leaving a ball of radius $\lceil t^\gamma \rceil$ for some $\gamma \in (0,1)$ but staying inside a ball of radius $t \log t$. We slice to annulus between these two balls into layers, and derive an estimate for paths that reach a given layer but do not reach the next layer.  To that end, fix $\gamma \in (\alpha,1)$ with $\alpha$ as in \eqref{e:def_Sr}, and let
\begin{equation}
\label{e:defrkt}
K_t := \lceil t^{1-\gamma} \log t \rceil, \qquad r^{(k)}_t := k \lceil t^\gamma \rceil, \quad
1 \leq k \leq K_t, \qquad \ell_t := K_t \lceil t^\gamma \rceil \geq t \log t.
\end{equation}
For $1 \leq k \leq K_t$, define (recall \eqref{defcurlyP})
\begin{equation}
\cN^{\ssup k}_t := \left\{ \pi \in \scrP(\cO, V) \colon\, \supp(\pi) \subset B_{r^{\ssup {k+1}}_t}(\cO),\, 
\supp(\pi)\cap B^\cc_{r^{\ssup k}_t}(\cO) \neq \emptyset \right\}
\end{equation}
and set
\begin{equation}
U^{\ssup k}(t) := \E_\cO \left[ \ee^{\int_0^t \xi(X_s) \dd s} \1_{\{\pi_{[0,t]}(X) \in \cN^{\ssup k}_t \}}\right].
\end{equation}

\begin{lemma}{\bf [Upper bound on $U^{\ssup k}(t)$]}
\label{l:UBpieces}
For any $\varepsilon>0$, $(\Prob \times \mathfrak{P})$-a.s.\ eventually as $t \to \infty$,
\begin{equation}
\label{e:UBpieces}
\sup_{1 \leq k \leq K_t} \frac 1t \log U^{\ssup k}_t \leq \frac{1}{t}\log U^*(t) + \varepsilon.
\end{equation}
\end{lemma}
\begin{proof}
We follow \cite[Lemma 4.4]{dHKdS2020}. Fix $k \in \{1, \ldots, K_t\}$. For $\pi \in \cN^{\ssup k}_t$, let
\begin{equation}
\gamma_\pi := \lambda_{r^{\ssup {k+1}}_t, A}(\pi) + \ee^{-S_{\lceil t^\gamma \rceil}}, 
\qquad z_\pi \in \supp(\pi), |z_\pi| > r^{\ssup k}_t,
\end{equation}
chosen such that \eqref{e:cond_massclass1}--\eqref{e:cond_massclass2} are satisfied. By Proposition~\ref{p:massclass} and \eqref{Aqdef},  $(\Prob \times \mathfrak{P})$-a.s.\ eventually as $t \to \infty$,
\begin{equation}
\label{e:prUBpieces2}
\begin{aligned}
\frac 1t \log U^{\ssup k}_t 
\leq \gamma_\pi - \frac{|z_\pi|}{t} \left( \log\log(\vartheta r^{(k+1)}_t)] 
- c  + o(1)  \right).
\end{aligned}
\end{equation}
Using Corollary~\ref{c:eigislandsGW} and $\log L_r \sim \vartheta r$, we bound
\begin{equation}
\begin{aligned}
\gamma_\pi 
\leq \varrho \log (\vartheta r^{(k+1)}_t) - \widetilde{\chi}(\varrho) + \tfrac12 \varepsilon + o(1).
\end{aligned}
\end{equation}
Moreover, $|z_\pi| > r^{\ssup {k+1}}_t - \lceil t^\gamma \rceil$ and
\begin{equation}
\begin{aligned}
&\frac{\lceil t^\gamma \rceil}{t} \left( \log\log( \vartheta r^{(k+1)}_t)] 
- c \right) \leq 
\frac{1}{t^{1-\gamma}} \log \log (2 t \log t) = o(1).
\end{aligned}
\end{equation}
Hence
\begin{equation}
\label{e:UBfinal}
\gamma_\pi \leq F_t(r^{(k+1)}_t) - \widetilde{\chi}(\varrho) + \tfrac12 \varepsilon + o(1)
\end{equation}
with
\begin{equation}
F_{c,t}(r) := \varrho \log (\vartheta r) - \frac{r}{t} \big[ \log\log (\vartheta r) - c \big], \qquad r>0.
\end{equation}
The function $F_{c,t}$ is maximised at any point $r_{c,t}$ satisfying
\begin{equation}
\label{keyasymp1}
\varrho t = r_{c,t}\log\log r_{c,t} - cr_{c,t} + \frac{r_{c,t}}{\log r_{c,t}}.
\end{equation}
In particular, $r_t = \mathfrak{r}_t[1+o(1)]$, which implies that
\begin{equation}
\label{e:prUBpieces1}
\sup_{r > 0} F_t(r) \leq \varrho \log (\vartheta \mathfrak{r}_t) - \varrho + o(1), \qquad t \to \infty.
\end{equation}
Inserting \eqref{e:prUBpieces1} into \eqref{e:UBfinal}, we obtain $\displaystyle \frac 1t \log U^{\ssup k}_t  < \varrho \log (\vartheta \mathfrak{r}_t) - \varrho - \widetilde{\chi}(\varrho) + \varepsilon$, which is the desired upper bound because $\varepsilon>0$ is arbitrary.
\end{proof}

\begin{proof}[Proof of the upper bound in \eqref{UU*comp}]
To avoid repetition, all statements hold $(\Probgr \times \Prob)$-a.s.\ eventually as $t \to \infty$. Set
\begin{equation}
U^{\ssup 0}(t) := \E_\cO \left[\ee^{\int_0^t \xi(X_s) \dd s} \1_{\{\tau_{[B_{\lceil t^\gamma \rceil}]^\cc}>t\}}\right],
\quad
U^{\ssup \infty}(t) := \E_\cO \left[ \ee^{\int_0^t \xi(X_s) \dd s} \1_{\{\tau_{[B_{\lceil t \log t \rceil}]^\cc} \leq t\}}\right].
\end{equation}
Then
\begin{equation}
U(t) \leq U^{\ssup 0}(t) + U^{\ssup \infty}(t) + K_t \max_{1 \leq k \leq K_t} U^{\ssup k}(t).
\end{equation}
From Lemmas~\ref{l:longpaths}--\ref{l:UBpieces} and the fact that $K_t = o(t)$, we get
\begin{equation}
\limsup_{t\to\infty} \left\{\frac{1}{t} \log U(t) - \frac{1}{t} \log U^*(t)\right\} \leq \varepsilon.
\end{equation}
Since $\varepsilon>0$ is arbitrary, this completes the proof of the upper bound in \eqref{e:QLyapGWT}.
\end{proof}


\subsection{Lower bound}
\label{sec:lb}
We follow \cite[Section 4.1]{dHKdS2020}. Fix $\varepsilon>0$. By the definition of $\widetilde{\chi}$, there exists an infinite rooted tree $T=(V',E',\YY)$ with degrees in $\supp(D_g)$ such that $\chi_T(\varrho) < \widetilde{\chi}(\varrho) + \tfrac14 \varepsilon$. Let $Q_r = B^T_r(\YY)$ be the ball of radius $r$ around $\YY$ in $T$. By Proposition~\ref{p:dualrepchi} and \eqref{e:defhatchi}, there exist a radius $R \in \N$ and a potential profile $q\colon B^T_R \to\R$ with $\cL_{Q_R}(q;\varrho)<1$ (in particular, $q\leq 0$) such that
\begin{equation}
\label{e:prLB0}
\lambda_{Q_R}(q;T) \geq -\widehat{\chi}_{Q_R}(\varrho;T) - \tfrac12 \varepsilon 
> -\widetilde{\chi}(\varrho) - \varepsilon.
\end{equation}
For $\ell\in\N$, let $B_\ell = B_\ell(\cO)$ denote the ball of radius $\ell$ around $\cO$ in $\GW$. We will show next that, $(\mathfrak{P} \times \Prob)$-a.s.\ eventually as $\ell \to \infty$, $B_\ell$ contains a copy of the ball $Q_R$ where the potentail $\xi$ is bounded from below by $\varrho\log\log |B_\ell(\cO)| + q$.

\begin{proposition}{\bf [Balls with high exceedances]}
\label{p:existencesubtree}
$(\Probgr \times \Prob)$-almost surely eventually as $\ell \to \infty$, there exists a vertex $z \in B_\ell$ with $B_{R+1}(z) \subset B_\ell$ and an isomorphism $\varphi:B_{R+1}(z) \to Q_{R+1}$ such that $\xi \geq \varrho \log \log |B_\ell(\cO)| + q \circ \varphi$ in $B_R(z)$. In particular,
\begin{equation}
\lambda_{B_R(z)}(\xi; \GW) > \varrho \log \log |B_\ell(\cO)| - \widetilde{\chi}(\varrho) - \varepsilon.
\end{equation}
Any such $z$ necessarily satisfies $|z| \geq c \ell$ $(\Probgr \times \Prob)$-a.s.\ eventually as $\ell \to \infty$ 
for some constant $c = c(\varrho, \vartheta, \widetilde{\chi}(\varrho), \varepsilon) >0$.
\end{proposition}

\begin{proof}
We follow \cite[Proposition 4.1]{dHKdS2020}. Only the last step changes as a result of the normalised Laplacian. 
First note that, as a consequence of the definition of $\GW$, it may be shown straightforwardly that, for some $p=p(T, R) \in (0,1)$ and $\Probgr$-almost surely eventually as $\ell \to \infty$, there exist $N \in \N$, $N \geq p |B_\ell|$ and distinct $z_1, \ldots, z_N \in B_\ell$ such that $B_{R+1}(z_i) \cap B_{R+1}(z_j) = \emptyset$ for $1\leq i\neq j\leq N$ and, for each $1 \leq i \leq N$, $B_{R+1}(z_i) \subset B_\ell$  and $B_{R+1}(z_i)$ is isomorphic to $Q_{R+1}$. Now, by \eqref{e:DE}, for each $i \in \{1,\ldots, N\}$,
\begin{equation}
\Prob\big(\xi\geq \varrho\log\log|B_\ell|+q \mbox{ in }B_R(z_i)\big)=|B_\ell|^{-\cL_{Q_R}(q)}.
\end{equation}
Using additionally that $|B_\ell| \geq \ell$ and $1-x \leq \ee^{-x}$, $x\in\R$, we obtain
\begin{equation*}
\begin{aligned}
\Prob(\not \exists i \in \{1,\ldots,N\} \colon \xi\geq \varrho\log\log|B_\ell|+q \mbox{ in }B_R(z_i))
= \left( 1 - |B_\ell|^{-\cL_{Q_R}(q)} \right)^N 
\leq \ee^{-p \ell^{1-\cL_{Q_R}(q)}},
\end{aligned}
\end{equation*}
which is summable in $\ell \in \N$, so the proof of the first statement is completed using the Borel-Cantelli lemma. As for the last statement note that by \eqref{e:monot_princev} and Lemma~\ref{l:maxpotential}
\begin{equation}
\lambda_{B_{c\ell}}(\xi;\GW) \leq \max_{x \in B_{c \ell}(\cO)} \xi(x) < a_{L_{c \ell}} +o(1) < a_{L_\ell} + \varrho \log c \vartheta +o(1) < a_{L_\ell} - \widetilde{\chi}(\varrho)-\varepsilon
\end{equation}
provided $c>0$ is small enough.
\end{proof}
\begin{lemma}
\label{lem:betterpath}
Let $z \in \GW$ and let $v_z = (v_{z,i})_{i=0}^{|z|}$  be the shortest path from $\cO$ to $z$, i.e. $v_{z,0} = \cO$, $v_{z,|z|} = z$, and $v_{z,i} \sim v_{z,i-1}$ for $i=1, \ldots, |z|$. Then
\begin{equation*}
\sum\limits_{L \in \mathbb{N}} \mathfrak{P}\left(\bigcup_{z \in Z_L} \left\{\prod_{i=1} ^L \frac{1}{\deg(v_{z,i})} \leq \frac{1}{(\log L)^{\delta_L L}}\right\} \right) < \infty,
\end{equation*}  
where $\delta_L$ satisfies $\lim\limits_{L\to \infty}\delta_L\log\log L = \infty$.
\end{lemma}
\begin{proof}
For $L \in \mathbb{N}$, let $\mathcal{Z}_L$ be the $L$-th generation of $\mathcal{GW}$ rooted at $\mathcal{O}$. For $z \in \mathcal{Z}_L$, let 
\begin{equation*}
\mathcal{E}_z = \left\{\prod_{i=1}^L \mathrm{deg}(v_{z_i}) \geq (\log L)^{\delta_L L}\right\}.
\end{equation*} 
We want to estimate
\begin{equation*}
\mathfrak{P}\left(\cup_{z \in \mathcal{Z}_L} \mathcal{E}_z\right). 
\end{equation*}
Pick any $K \in \mathbb{N}$ and estimate
\begin{equation*}
\mathfrak{P}\left(\cup_{z \in \mathcal{Z}_L} \mathcal{E}_z\right)
\leq \mathfrak{P}\left(\cup_{z \in \mathcal{Z}_L} \mathcal{E}_z,\,|\mathcal{Z}_L| > K\right)
+ \mathfrak{P}\left(\cup_{z \in \mathcal{Z}_L} \mathcal{E}_z,\,|\mathcal{Z}_L| \leq K\right).
\end{equation*}
Estimate
\begin{equation*}
\mathfrak{P}\left(\cup_{z \in \mathcal{Z}_L} \mathcal{E}_z,\,|\mathcal{Z}_L| > K\right)
\leq \mathfrak{P}\left(|\mathcal{Z}_L| > K\right) \leq \frac{1}{K} \mathfrak{E}(\mathcal{Z}_L) =: \frac{e_L}{K}.
\end{equation*}
Also estimate
\begin{equation*}
\begin{aligned}
&\mathfrak{P}\left(\cup_{z \in \mathcal{Z}_L} \mathcal{E}_z,\,|\mathcal{Z}_L| \leq K\right)\\
&= \sum_{\ell=1}^K \mathfrak{P}\left(\cup_{k=1}^\ell \mathcal{E}_{z_\ell},\,|\mathcal{Z}_L| = \ell\right)
\leq \sum_{\ell=1}^K \sum_{k=1}^\ell \mathfrak{P}\left(\mathcal{E}_{z_\ell},\,|\mathcal{Z}_L| = \ell\right)\\
&\leq \sum_{\ell=1}^K \sum_{k=1}^\ell \mathfrak{P}\left(\mathcal{E}_{z_1},\,|\mathcal{Z}_L| = \ell\right)
\leq K \sum_{\ell=1}^K \mathfrak{P}\left(\mathcal{E}_{z_1},\,|\mathcal{Z}_L| = \ell\right)\\
&= K \mathfrak{P}\left(\mathcal{E}_{z_1},\,|\mathcal{Z}_L| \leq K\right) 
\leq K \mathfrak{P}\left(\mathcal{E}_{z_1}\right) =: Kp_L,
\end{aligned}
\end{equation*}
where $z_\ell$ is the $\ell$-th vertex in $\mathcal{Z}_L$ (say in lexicographic order), and $p_L$ is the probability that the product of $L$ i.i.d.\ copies of the degrees exceeds $(\log L)^{\delta_L L}$. In the last inequality we need not worry about the correlation between $\mathcal{E}_{z_1}$ and the event $|\mathcal{Z}_L| \leq K$ because we drop the latter. Thus, for any $K \in \mathbb{N}$ we have
\begin{equation*}
\mathfrak{P}\left(\cup_{z \in \mathcal{Z}_L} \mathcal{E}_z\right) \leq \frac{e_L}{K} + Kp_L.
\end{equation*}
Now minimise over $K$. The minimising value is $K= \sqrt{e_L/p_L}$ (to be rounded off to an integer), so that we get   
\begin{equation*}
\mathfrak{P}\left(\cup_{z \in \mathcal{Z}_L} \mathcal{E}_z\right) \leq 2 \sqrt{e_Lp_L}.
\end{equation*}
Since $e_L = \mathrm{e}^{\vartheta L + o(L)}$ and $p_L = \mathrm{e}^{-L\delta_L \log\log L + O(L)}$, it follows that 
\begin{equation*}
\sum_{L \in \mathbb{N}} \mathfrak{P}\left(\cup_{z \in \mathcal{Z}_L} \mathcal{E}_z\right) < \infty
\end{equation*} 
by the assumption on $\delta_L$.
\end{proof}
\noindent Lemma~\ref{lem:betterpath} implies that $\mathfrak{P}$-almost surely eventually as $L \to \infty$, any path $y_z$ must satisfy
\begin{equation}
\prod_{i=1} ^L \frac{1}{\deg(y_{z,i})} \geq \frac{1}{(\log L)^{\delta_L L}}.  
\end{equation}
\begin{proof}[Proof of the lower bound in \eqref{e:QLyapGWT}]
Let $z$ be as in Proposition~\ref{p:existencesubtree}. Write $\tau_z$ for the hitting time of $z$ by the random walk $X$. For $s\in (0,t)$, we estimate
\begin{equation}
\label{lowbound1}
\begin{aligned}
U(t) &\geq \E_\cO\Big[\ee^{\int_0^t \xi(X_u)\,\dd u}\,\1_{\{\tau_z\leq s\}}\,
\1_{\{X_u\in B_R(z)\,\forall u\in[\tau_z,t]\}}\Big]\\
&=\E_\cO\Big[\ee^{\int_0^{\tau_z} \xi(X_u)\,\dd u}\,\1_{{\{\tau_z\leq s\}}}\,
\E_z\Big[\ee^{\int_0^{v} \xi(X_u)\,\dd u}\,\1_{{\{X_u\in B_R(z)\,\forall u\in [0,v]\}}}\Big]\Big|_{v=t-\tau_z}\Big],
\end{aligned}
\end{equation}
where we use the strong Markov property at time $\tau_z$. We first bound the last term in the integrand in \eqref{lowbound1}. Since $\xi \geq \varrho \log\log |B_\ell| +q $ in $B_R(z)$, 
\begin{equation}
\begin{aligned}
\E_z\Big[\ee^{\int_0^{v} \xi(X_u)\,\dd u} \1_{\{X_u\in B_R(z)\,\forall u\in [0,v]\}}\Big]
& \geq \ee^{v \varrho \log \log |B_\ell|} \E_{\YY}\Big[\ee^{\int_0^{v} q(X_u)\,\dd u} 
\1_{\{X_u\in Q_R\,\forall u\in [0,v]\}}\Big] \\
& \geq \e^{ v \varrho \log \log |B_\ell|} \ee^{v \lambda_{Q_R}(q;T)} \phi^{\ssup 1}_{Q_R}(\YY)^2 \\
& > \exp \big\{ v \left(\varrho \log\log |B_{\ell}| -  \widetilde{\chi}(\varrho) - \varepsilon \right) \big\}
\end{aligned}
\end{equation}
for large $v$, where we use that $B_{R+1}(z)$ is isomorphic to $Q_{R+1}$ for the indicators in the first inequality, and apply Lemma~\ref{l:bounds_mass} and \eqref{e:prLB0} to obtain the second and third inequalities respectively. On the other hand since $\xi\geq0$, we have  
\begin{equation}
\E_\cO \Big[\ee^{\int_0^{\tau_z} \xi(X_u)\,\dd u}\1{\{\tau_z\leq s\}}\Big]
\geq \P_\cO(\tau_z\leq s),
\end{equation}
and we can bound the latter probability from below by the probability that the random walk runs along a shortest path from the root $\cO$ to $z$ within a time at most $s$. This gives
\begin{equation}
\begin{aligned}
\Prob_\cO(\tau_z\leq s) 
& \geq \Big(\prod_{i=1}^{|z|}\frac 1{\deg(y_{z,i})}\Big) P\Big(\sum_{i=1}^{|z|} E_i \leq s\Big) 
\geq (\log |z|)^{-\delta_{|z|}|z|}
{\rm Poi}_{d_{\rm min} s}([|z|,\infty)),
\end{aligned}
\end{equation}
where ${\rm Poi}_\gamma$ is the Poisson distribution with parameter $\gamma$, and $P$ is the generic symbol for probability. The final inequality uses Lemma~\ref{lem:betterpath}. Summarising, we obtain
\begin{equation}
\label{lowbound2}
\begin{aligned}
U(t) 
& \geq (\log |z|)^{-\delta_{|z|}|z|} \e^{-d_{\rm min} s}\frac{(d_{\rm min} s)^{|z|}}{|z|!}
\e^{(t-s)\left[\varrho\log\log |B_{\ell}| - \widetilde{\chi}(\varrho) - \varepsilon \right]} \\
& \geq \exp \left\{-d_{\min} s + (t-s)\left[\varrho\log\log |B_{\ell}| - \widetilde{\chi}(\varrho) 
- \varepsilon \right] - |z| \log \left( \frac{(\log |z|)^{\delta_{|z|}}}{d_{\min}}\frac{|z|}{s}\right) \right\} \\
& \geq \exp \left\{-d_{\min} s + (t-s)\left[\varrho\log\log |B_{\ell}| - \widetilde{\chi}(\varrho) 
- \varepsilon \right] - \ell \log \left( \frac{(\log \ell)^{\delta_{\ell}}}{d_{\min}}\frac{\ell}{s}\right) \right\},
\end{aligned}
\end{equation}
where in the last inequality we use that $s \leq |z|$ and $\ell \geq |z|$. Further assuming that $\ell = o(t)$, we see that the optimum over $s$ is obtained at 
\begin{equation}
s= \frac{\ell}{d_{\min}+\varrho\log\log|B_{\ell}|-\widetilde{\chi}(\varrho) - \varepsilon} =o(t).
\end{equation} 
Note that, by Proposition~\ref{p:existencesubtree}, this $s$ indeed satisfies $s\leq |z|$. Applying \eqref{e:volumerateGW} we get, after a straightforward computation, $(\Probgr \times \Prob)$-a.s.\ eventually as $t \to \infty$,
\begin{equation}
\label{Ulowboundwithr}
\frac 1t\log U(t) \geq \varrho\log \log |B_\ell| - \frac{\ell}{t} \log\log \ell -  \frac{\ell}t \delta_\ell \log \log \ell 
- \widetilde{\chi}(\varrho) - \varepsilon + O\left( \frac{\ell}{t} \right).
\end{equation}
Inserting $\log |B_\ell| \sim \vartheta \ell$, we get
\begin{equation}
\frac 1t\log U(t) \geq F_\ell - \widetilde{\chi}(\varrho) - \varepsilon + o(1) + O\left( \frac{\ell}{t} \right) 
\end{equation}
with
\begin{equation}
F_\ell = \varrho\log(\vartheta \ell) - \frac{\ell}{t} \log\log \ell -  \frac{\ell}t \delta_\ell \log \log \ell. 
\end{equation}
The optimal $\ell$ for $F_\ell$ satisfies
\begin{equation}
\label{keyasymp2}
\varrho t = \ell[1+\delta_\ell + l\tfrac{\dd}{\dd \ell} \delta_\ell ] \log\log\ell 
+ \frac{\ell\delta_\ell}{\log \ell}  + \frac{\ell}{\log \ell}, 
\end{equation}
i.e. $\ell = \mathfrak{r}_t[1+o(1)]$. For this choice we obtain
\begin{equation}
\label{UlowboundGWT1}
\frac 1t\log U(t)\geq \varrho\log(\vartheta\mathfrak{r}_t) - \varrho
-\widetilde{\chi}(\varrho) - \varepsilon + o(1).
\end{equation}
Hence $(\Probgr \times \Prob)$-a.s.
\begin{equation}
\label{UlowboundGWT2}
\liminf_{t \to \infty}
\left\{ \frac{1}{t}\log U(t) - \frac{1}{t} \log U^*(t)\right\} \geq - \varepsilon.
\end{equation}
Since $\varepsilon>0$ is arbitrary, this completes the proof of the lower bound in \eqref{e:QLyapGWT}.
\end{proof}

\section{Analysis of the variational formula}
\label{s:analysischi}
This section is dedicated the analysis of variational formula. Proposition~\ref{p:dualrepchi} is proven in Section~\ref{proofsection1}. Theorem~\ref{t:tildechilargerho} is proven in Section~\ref{ss:theoremvarproof}, which is done by adapting the gluing argument in \cite{dHKdS2020}.
\subsection{Proof of Proposition~\ref{p:dualrepchi}}
\label{proofsection1}
The inequality is clear. For the equality we first prove that for any graph $G = (V,E)$ and $\Lambda \subset V$,
\begin{equation}
\label{e:hatchidualrep1}
\widehat{\chi}_\Lambda(\varrho; G) = \inf_{\substack{p \in \cP(V)\colon\, \\ \supp(p) \subset \Lambda} } 
\left[ I_E(p) + \varrho J_V(p) \right].
\end{equation}
For this we follow \cite[Lemma 2.17]{GM1998}.
By the Rayleigh-Ritz formula,
\begin{align*}
    \lambda_{\Lambda}(q;G)&=\sup_{{\supp(\phi)\subset \Lambda} \atop {\lVert \phi \rVert=1 }} \langle (\Delta_{G}+q)\phi , \phi \rangle =
    \sup_{\lVert \phi \rVert=1 } \left\{\sum_{x\in \Lambda}  \deg(x)[(\Delta\phi)(x)+q(x)\phi(x)] \phi(x) \right\}\\
    &=\sup_{\lVert \phi \rVert=1 } \left\{-\sum_{\{x, y\}\in E_{\Lambda}}[\phi(x)-\phi(y)]^{2} + \sum_{x\in \Lambda} 
 \deg(x) q(x)\phi(x)^{2} \right\}.
\end{align*}
By Lemma~\ref{lem:specbd}(2), the eigenfunction corresponding to $\lambda_\Lambda(q;G)$ may be taken to be non-negative, and we may therefore make the substitution $\phi(x)=\sqrt{\frac{p(x)}{\deg(x)}}$ so that  $p$ with $p(x) = \phi(x)^2\deg(x)$ is a probability measure supported on $\Lambda$. So
\begin{align*}
    \lambda_{\Lambda}(q;G)=\sup_{{p\in \mathcal{P}(V)} \atop  {\supp(p)\subset \Lambda }} \left\{-I_{E_{\Lambda}}(p) +\sum_{x\in\Lambda } q(x) p(x)\right\},
\end{align*}
and therefore
\begin{align*}
    \widehat{\chi}_{\Lambda}(\varrho; G) &=
    -\sup_{{q:V\to [-\infty,\infty) } \atop {\mathcal{L}_{V}(q;\varrho)}} \left[ \sup_{p\in \mathcal{P}(\Lambda)} \left\{-I_{E_\Lambda}(p)+\sum_{x\in \Lambda} q(x)p(x) \right\} \right]\\
    &= -\sup_{p\in \mathcal{P}(\Lambda)}\left[ \sup_{q\colon \mathcal{L}(q;\varrho)=1} \left\{ \sum_{x\in \Lambda}q(x)p(x)-\varrho \log \sum_{x\in \Lambda}\ee^{q(x)/\varrho}  \right\}- I_{E_\Lambda}(p) \right].
\end{align*}
As the expression in the curly brackets does not change by adding a constant to $q(x)$, the inner supremum may be taken over all $q: \Lambda \to \mathbb{R}$. 

For $z\in \Lambda$, differentiating with respect to $q(z)$ and setting equal to $0$ we get that the supremum is attained at $\Bar{q}$ satisfying
\begin{equation*}
p(z) = \frac{\ee^{\Bar{q}(z)/\varrho}}{\sum_{x \in \Lambda} \ee^{\Bar{q}(x)/\varrho}}  
\end{equation*}
for all $z$. Or equivalently,
\begin{equation*}
\Bar{q}(z) = \varrho \log p(z) + \varrho \log \sum_{x \in \Lambda} \ee^{\Bar{q}(x)/\varrho}. 
\end{equation*}
This gives that the value of the inner supremum is $-\varrho J_V(p)$ and \eqref{e:hatchidualrep1} follows.\\

Recall the definition of $B_r(\cO)$ from \eqref{e:volumerateGW}. By \eqref{e:hatchidualrep1}, $\widehat{\chi}_{B_r(\cO)}(\varrho; G)$ is non-increasing in $r$ and therefore,
\begin{equation}
\label{e:chi0small}
\lim_{r \to \infty} \widehat{\chi}_{B_r(\cO)}(\varrho ; G) \geq \chi_G(\varrho).
\end{equation}    
It remains to show the opposite inequality. For that we show that for any $p\in \cP(V)$ and $r \in \N$, there exists a $p_r\in\cP(V)$ with support in $B_r(\cO)$ such that 
\begin{equation}
\label{e:prchi0large1}
\liminf_{r\to\infty} \left\{ I_E(p_r)+\varrho J_V(p_r) \right\} \leq I_E(p)+\varrho J_V(p).
\end{equation} 
We follow \cite[Lemma A.2]{dHKdS2020}. Simply take
\begin{equation}
p_r(x)= \frac{p(x)\1_{B_r(\cO)}(x)}{p(B_r(\cO))}, \qquad x \in V,
\end{equation} 
i.e. the normalized restriction of $p$ to $B_r(\cO)$. Then we easily see that 
\begin{equation}
\begin{aligned}
J_V(p_r)-J_V(p) &=-\frac1{p(B_r(\cO))} \sum_{x\in B_r(\cO)} p(x)\log p(x)+\log p(B_r(\cO))+\sum_{x\in V} p(x)\log p(x)\\
& \leq  \frac {J_V(p)}{p(B_r(\cO))}(1-p(B_r(\cO))) \underset{r \to \infty}{\longrightarrow} 0,
\end{aligned}
\end{equation}
where we use that $\log p(B_r(\cO)) \leq 0$ and $p(x) \log p(x)\leq 0$ for every $x$. As for the $I$-term,
\begin{equation}
\begin{aligned}
I_E(p_r)&=\frac1{p(B_r(\cO))}\sum_{\{x,y\}\in E \colon x,y\in B_r(\cO)}\left(\sqrt{\tfrac{p(x)}{\deg(x)}}-\sqrt{\tfrac{p(y)}{\deg(y)}}\,\right)^2\\
&\quad +\frac1{p(B_r(\cO))}\sum_{\{x, y\} \in E \colon x\in B_r(\cO),\, y\in B^\cc_r}\frac{p(x)}{\deg(x)}
\leq \frac{I_E(p)}{p(B_r(\cO))}+\frac{p(B_{r-1}(\cO)^\cc)}{d_{\min}\ p(B_r(\cO))},
\end{aligned}
\end{equation}
and therefore
\begin{equation}
I_E(p_r)-I_E(p) \leq \frac{I_E(p)}{p(B_r(\cO))}(1-p(B_r(\cO)))+\frac{p(B_{r-1}(\cO)^\cc)}{d_{\min}\ p(B_r(\cO))} \underset{r \to \infty}{\longrightarrow} 0.
\end{equation}

\subsection{Proof of Theorem~\ref{t:tildechilargerho}}
\label{ss:theoremvarproof}
This section follows \cite[Appendix A]{dHKdS2020}. We adapt the techniques to the new $I_E$ function defined in \eqref{e:defIJ}. 
\begin{lemma}{\bf [Glue two]}
\label{l:compjoin2}
Let $G_i = (V_i, E_i)$, $i=1,2$, be two disjoint connected simple graphs, and let $x_i \in V_i$, $i=1,2$. Denote by ${G}$
the union graph of $G_1$, $G_2$ with one extra edge between $x_1$ and $x_2$, i.e. $ {G} = ( {V},  {E})$ with $ {V} := V_1 \cup V_2$, $ {E} := E_1 \cup E_2 \cup \{(x_1, x_2)\}$. Then
\begin{equation}
\label{e:glueingchi}
\chi_{ {G}} \geq \min \left\{ \chi_{G_1} , \chi_{G_2}  \right\}.
\end{equation}
\end{lemma}

\begin{proof} 
We follow \cite[Lemma A.3]{dHKdS2020}. Given $p \in \cP(V)$, let $a_i = p(V_i)$, $i=1,2$, and define $p_i \in \cP(V_i)$ by putting
\begin{equation}
p_i(x) := \begin{cases}
\frac{1}{a_i} p(x) \1_{V_i}(x) & \text{if } a_i> 0,\\
\1_{x_i}(x) & \text{otherwise.}
\end{cases}
\end{equation} 
Straightforward manipulations show that
\begin{equation}
I_E(p) = \sum_{i=1}^2 a_i I_{E_i}(p_i) + \left(\sqrt{\tfrac{p(x_1)}{\deg(x_1)}} - \sqrt{\tfrac{p(x_2)}{\deg(x_2)}}\right)^2, 
\qquad J_V(p) = \sum_{i=1}^2 \left[ a_i J_{V_i}(p_i) - a_i \log a_i \right],
\end{equation}
and so
\begin{equation}
I_E(p) + \varrho J_V(p) \geq \sum_{i=1}^2 a_i \Big[ I_{E_i}(p_i) 
+ \varrho J_{V_i}(p_i) \Big]  \geq \min\{\chi_{G_1}, \chi_{G_2}\}.
\end{equation}
The proof is completed by taking the infimum over $p \in \cP(V)$.
\end{proof}
Below it will be useful to define, for $x \in V$,
\begin{equation}
\label{e:altchirest}
\chi^{\ssup{x,b}}_G  = \inf_{\substack{p \in \cP(V),\\p(x)=b}} [I_E(p) + \varrho J_V(p)],
\end{equation}
i.e. a version of $\chi_G$ with ``boundary condition'' $b$ at $x$. It is clear that $\chi^{\ssup{x,b}}_G \geq \chi_G$.
Next we glue several graphs together and derive representations and estimates for the corresponding $\chi$. For $k \in \N$, let $G_i = (V_i, E_i)$, $1 \leq i \leq k$, be a collection of disjoint graphs. Let $x$ be a point not belonging to $\bigcup_{i=1}^k V_i$. For a fixed choice $y_i \in V_i$, $1 \leq i \leq k$, we denote by $\overline{G}_k = (\overline{V}_k, \overline{E}_k)$ the graph obtained by adding an edge from each $y_1, \ldots, y_k$ to $x$, i.e. $\overline{V}_k = V_1 \cup \cdots \cup V_k \cup \{x\}$ and $\overline{E}_k = E_1 \cup \cdots \cup E_k \cup  \{(y_1, \cO),\dots,(y_k,x)\}$.

\begin{lemma}{\bf [Glue many plus vertex]}
\label{l:joiningk}
For any $\varrho>0$, any $k \in \N$, and any $G_i=(V_i, E_i)$, $y_i \in V_i$, $1\leq i \leq k$,
\begin{equation}
\label{e:joiningk}
\begin{aligned}
\chi_{\overline{G}_{k}}  = 
& \inf_{\substack{0 \leq c_i \leq a_i \leq 1, \\ a_1 + \cdots + a_k \leq 1}}
\Big\{
 \sum_{i=1}^k a_i \left(\chi_{G_i}^{\ssup{y_i, c_i/a_i}} - \varrho \log a_i \right) \\
& \;\; + \sum_{i=1}^k \left( \sqrt{\frac{c_i}{\deg(y_i)}} - \sqrt{\frac{1-\sum_{i=1}^k a_i}{\deg(x)}} \right)^2
 - \varrho \Big(1-\sum_{i=1}^k a_i \Big) \log \Big(1-\sum_{i=1}^k a_i\Big)
\Big\}.
\end{aligned}
\end{equation}
\end{lemma}
\begin{proof}
We follow \cite[Lemma A.4]{dHKdS2020}.
The claim follows from straightforward manipulations with \eqref{e:defIJ}.
\end{proof}
\noindent Lemma~\ref{l:joiningk} leads to the following comparison lemma. For $j \in \N$, let
\begin{equation}
(G^j_i,y^j_i) = \begin{cases}
(G_i, y_i) & \text{if } i < j,\\
(G_{i+1}, y_{i+1}) & \text{if } i \geq j,
\end{cases}
\end{equation}
i.e. $(G^j_i)_{i\in\N}$ is the sequence $(G_i)_{i\in\N}$ with the $j$-th graph omitted. Let $\overline{G}^j_k$ be the analogue of $\overline{G}_k$ obtained from $G^j_i$, $1 \leq i \leq k$, $i\neq j$, instead of $G_i$, $1 \leq i \leq k$.
\begin{lemma}{\bf [Comparison]}
\label{l:comparejoinedk}
For any $\varrho>0$ and any $k \in \N$,
\begin{equation}
\label{e:joinedk+1}
\begin{aligned}
\chi_{\overline{G}_{k+1}}  
& = \inf_{1 \leq j \leq k+1} \; \inf_{0 \leq c \leq u \leq \tfrac{1}{k+1}} \; 
\inf_{\substack{0 \leq c_i \leq a_i \leq 1, \\ a_1+ \cdots +  a_k \leq 1}}
\Bigg\{ 
(1-u) \Big[
\sum_{i=1}^k a_i \big(\chi_{G_{\sigma_j(i)}}^{\ssup{y_{\sigma_j(i)}, c_i/a_i}} - \varrho \log a_i \big) \\
&  \qquad + \sum_{i=1}^k\left( \sqrt{\frac{c_i}{\deg(y_i)}} - \sqrt{\frac{1-\sum_{i=1}^k a_i}{\deg(x)}} \right)^2 
 - \varrho \Big(1-\sum_{i=1}^k a_i \Big) \log \Big(1-\sum_{i=1}^k a_i \Big) \Big]  \\
&  \qquad + u \chi_{G_j}^{\ssup{y_j, c/u}} + \left(\sqrt{\frac{c}{\deg(y_j)}} - \sqrt{\frac{(1-u)\Big(1-\sum_{i=1}^k a_i \Big)}{\deg(x)} } \right)^2 \\
& \qquad - \varrho \left[ u \log u + (1-u) \log (1-u) \right] 
\Bigg\}.
\end{aligned}
\end{equation}
Moreover, 
\begin{equation}
\label{e:comparedjoinedk}
\begin{aligned}
\chi_{\overline{G}_{k+1}} 
& \geq \inf_{1\leq j \leq k+1} \inf_{ 0 \leq u \leq \tfrac{1}{k+1}} \Bigg\{ (1-u) \chi_{\overline{G}^j_k}  \\
& \qquad \qquad \qquad \qquad \quad  + \inf_{v \in [0,1]} \Big\{ u \chi_{G_j}^{(y_j, v)} + \mathbbm{1}_{\{u(1+v) \geq 1 \}}\Big[\sqrt{\tfrac{vu}{\deg(y_j)}} - \sqrt{\tfrac{1-u}{\deg(x)}}\, \Big]^2 \Big\} \\
& \qquad \qquad \qquad \qquad \quad  - \varrho \left[u \log u + (1-u) \log(1-u) \right]
\Bigg\}.
\end{aligned}
\end{equation}
\end{lemma}
\begin{proof}
See \cite[Lemma A.5]{dHKdS2020}. The argument still applies with the definition of $I_E$ given in \eqref{e:defIJ}.
\end{proof}
\begin{lemma}{\bf [Propagation of lower bounds]}
\label{l:comparelargerho}
If $\varrho>0$, $M \in \R$,  $C >0$ and $k \in \N$ satisfy $\varrho \geq C/\log(k+1)$ and
\begin{equation}
\label{e:assumpcomparelargerho}
\inf_{1 \leq j \leq k+1} \chi_{\overline{G}^j_k}  \geq M, \qquad 
\inf_{1 \leq j \leq k+1} \inf_{v \in [0,1]} \chi_{G_j}^{\ssup{y_j, v}}  \geq M-C,
\end{equation}
then $\chi_{\overline{G}_{k+1}}  \geq M$.
\end{lemma}
\begin{proof}
See \cite[Lemma A.6]{dHKdS2020}. The proof carries over directly since $I_E$ does not appear.
\end{proof}
\noindent
The above results will be applied in the next section to minimise $\chi$ over families of trees with minimum degrees.
\subsubsection{Trees with minimum degrees}
\label{ss:treesmindeg}

Fix $d \in \N$. Let $\mr{\cT}_d$ be an infinite tree rooted at $\cO$ such that the degree of $\cO$ equals $d-1$ and the degree of every other vertex in $\mr{\cT}_d$ is $d$. Let $\mathring{\mathscr{T}}_d^{\ssup 0} = \{\mathring{\cT}_d\}$ and, recursively, let $\mr{\mathscr{T}}_d^{\ssup{n+1}}$ denote the set of  all trees obtained from a tree in $\mr{\mathscr{T}}_d^{\ssup n}$ and a disjoint copy of $\mr{\cT}_d$ by adding an edge between a vertex of the former and the root of the latter. Write $\mr{\mathscr{T}}_d = \bigcup_{n \in\N_0} \mr{\mathscr{T}}_d^{\ssup n}$. Assume that all trees in $\mr{\mathscr{T}}_d$ are rooted at $\cO$.

Recall that $\cT_d$ is the infinite regular $d$-tree. Observe that $\cT_d$ is obtained from $(\mr{\cT}_d, \cO)$ and a disjoint copy $(\mr{\cT}_d', \cO')$ by adding one edge between $\cO$ and $\cO'$. Consider $\cT_d$ to be rooted at $\cO$. Let $\mathscr{T}_d^{\ssup 0} = \{\cT_d\}$ and, recursively, let $\mathscr{T}_d^{\ssup{n+1}}$ denote the set of all trees obtained from a tree in $\mathscr{T}_d^{\ssup n}$ and a disjoint copy of $\mr{\cT}_d$ by adding an edge between a vertex of the former and the root of the latter. Write $\mathscr{T}_d = \bigcup_{n \in\N_0} \mathscr{T}_d^{\ssup n}$, and still consider all trees in $\mathscr{T}_d$ to be rooted at $\cO$. Note that $\mathscr{T}_d^{\ssup n}$ contains precisely those trees of $\mr{\mathscr{T}}_d^{\ssup{n+1}}$ that have $\cT_d$ as a subgraph rooted at $\cO$. In particular, $\mathscr{T}_d^{\ssup n} \subset \mr{\mathscr{T}}_d^{\ssup{n+1}}$ and $\mathscr{T}_d \subset \mr{\mathscr{T}}_d$.

Our objective is to prove the following.

\begin{proposition}{\bf [Minimal tree is optimal]}
\label{p:minimalchiT}
If $\varrho \geq \tfrac{1}{(d-1)\log(d+1)}$, then 
\[
\chi_{\cT_d}(\varrho)  = \min_{T \in \mathscr{T}_d} \chi_T(\varrho).
\]
\end{proposition}

For the proof of Proposition~\ref{p:minimalchiT}, we will need the following. 

\begin{lemma}{\bf [Minimal half-tree is optimal]}
\label{l:minimalchimrT}
For all $\varrho \in (0,\infty)$, 
\[
\chi_{\mr{\cT}_d}(\varrho)  = \min_{T \in \mr{\mathscr{T}}_d} \chi_T(\varrho).
\]
\end{lemma}

\begin{proof}
See \cite[Lemma A.8]{dHKdS2020}. The proof carries over directly since $I_E$ does not appear.
\end{proof}
\begin{lemma}{\bf [A priori bounds]}
\label{l:condchiT}
For any $d \in \N$ and any $\varrho \in (0,\infty)$,
\begin{equation}
\chi_{\mr{\cT}_d}(\varrho)  \leq \chi_{\cT_d}(\varrho)  \leq \chi_{\mr{\cT}_d}(\varrho)  + \frac{1}{d-1}.
\end{equation}
\end{lemma}

\begin{proof}
We follow \cite[Lemma A.9]{dHKdS2020}. The first inequality follows from Lemma~\ref{l:minimalchimrT}. For the second inequality, note that $\cT_d$ contains as subgraph a copy of $\mr{\cT}_d$, and restrict the minimum in \eqref{e:defchiG} to $p \in \cP(\mr{\cT}_d)$. 
\end{proof}
\begin{proof}[Proof of Proposition~\ref{p:minimalchiT}]
We follow \cite[Proposition A.7]{dHKdS2020}. Fix $\varrho \geq \tfrac{1}{(d-1)\log(d+1)}$. It will be enough to show that
\begin{equation}
\label{e:prminchiT1}
\chi_{\cT_d}  = \min_{T \in \mathscr{T}_d^{\ssup n}}\chi_T,\qquad n\in\N_0 .
\end{equation}
We will prove this by induction in $n$. The case $n=0$ is trivial. Assume that, for some $n_0\geq 0$, \eqref{e:prminchiT1} holds for all $n \leq n_0$. Let $T \in \mathscr{T}^{\ssup{n_0+1}}_d$. Then there exists a vertex $x$ of $T$ with degree $k+1 \geq d+1$. Let $y_1, \ldots, y_{k+1}$ be set of neighbours of $x$ in $T$. When we remove the edge between $y_j$ and $x$, we obtain two connected trees; call $G_j$ the one containing $y_j$, and $\overline{G}^j_k$ the other one. With this notation, $T$ may be identified with $\overline{G}_{k+1}$. 

Now, for each $j$, the rooted tree $(G_j,y_j)$ is isomorphic (in the obvious sense) to a tree in $\mr{\mathscr{T}}_d^{\ssup{\ell_j}}$, where $\ell_j \in \N_0$ satisfy $\ell_1 + \cdots + \ell_{k+1} \leq n_0$, while $\overline{G}^j_k$ belongs to $\mathscr{T}_d^{(n_j)}$ for some $n_j \leq n_0$. Therefore, by the induction hypothesis,
\begin{equation}
\label{e:prminchiT2}
\chi_{\overline{G}^j_k}  \geq \chi_{\cT_d},
\end{equation}
while, by \eqref{e:altchirest}, Lemma~\ref{l:minimalchimrT} and Lemma~\ref{l:condchiT},
\begin{equation}
\label{e:prminchiT3}
\inf_{v \in [0,1]} \chi_{G_j}^{(y_j, v)}  \geq \chi_{G_j}  \geq \chi_{\mr{\cT}_d}  \geq \chi_{\cT_d}  -\frac{1}{d-1}.
\end{equation}
Thus, by Lemma~\ref{l:compjoin2} applied with $M=\chi_{\cT_d} $ and $C=\tfrac{1}{d-1}$,
\begin{equation}
\chi_{T}  = \chi_{\bar{G}_{k+1}}  \geq \chi_{\cT_d}, 
\end{equation}
which completes the induction step.
\end{proof}
\begin{proof}[Proof of Theorem~\ref{t:tildechilargerho}]We follow \cite[Theorem 1.2]{dHKdS2020}.
First note that, since $\cT_{d_{\min}}$ has degrees in $\supp(D_g)$, $\widetilde{\chi}(\varrho) \leq \chi_{\cT_{d_{\min}}}(\varrho)$. For the opposite inequality, we proceed as follows. Fix an infinite tree $T$ with degrees in $\supp(D_g)$, and root it at a vertex $\YY$. For $r \in \N$, let $\widetilde{T}_r$ be the tree obtained from $B_r(\cO)= B^T_r(\YY)$ by attaching to each vertex $x \in B_r(\cO)$ with $|x| = r$ a number $d_{\min}-1$ of disjoint copies of $(\mathring{\cT}_{d_{\min}}, \cO)$, i.e. adding edges between $x$ and the corresponding roots. Then $\widetilde{T}_r \in \mathscr{T}_{d_{\min}}$ and, since $B_r(\cO)$ has more out-going edges in $T$ than in $\widetilde{T}_r$, we may check using \eqref{e:hatchidualrep1} that
\begin{equation}
\widehat{\chi}_{B_r}(\varrho ; T) \geq \widehat{\chi}_{B_r}(\varrho; \widetilde{T}_r) 
\geq \chi_{\widetilde{T}_r}(\varrho) \geq \chi_{\cT_{d_{\min}}}(\varrho).
\end{equation}
Taking $r \to \infty$ and applying Proposition~\ref{p:dualrepchi}, we obtain $\chi_T(\varrho) \geq \chi_{\cT_{d_{\min}}}(\varrho)$. Since $T$ is arbitrary, the proof is complete.
\end{proof}


\appendix





\begin{thebibliography}{99}

\bibitem{dHKdS2020}
F.\ den Hollander, W.\ Konig, R.S.\ dos Santos,
The Parabolic Anderson model on a Galton-Watson tree,
to appear in \emph{In and Out of Equilibrium 3: Celebrating Vladas Sidoravicius}, Progress in Probability, Birkh\"auser, 2021.

\bibitem{AGH2020}
L.\ Avena, O.\ G{\"u}n, M.\ Hesse, 
The parabolic Anderson model on the hypercube, 
Stoch.\ Proc.\ Appl.\ 130, 3369--3393, 2020.


\bibitem{BK2016}
M.\ Biskup, W.\ K\"onig,
Eigenvalue order statistics from random Schr\"odinger operators with doubly-exponential tails,
Commun.\ Math.\ Phys.\ 341, 179--218, 2016.

\bibitem{FM1990}
K.\ Fleischmann, S.A.\ Molchanov,
Exact asymptotics in a mean field model with random potential,
Probab.\ Theory Relat.\ Fields 86, 239--251, 1990. 

\bibitem{GM1990}
J.\ G\"artner, S.A.\ Molchanov,
Parabolic problems for the Anderson model I. Intermittency and related problems,
Commun.\ Math.\ Phys.\ 132,  613--655, 1990.

\bibitem{GM1998}
J.\ G\"artner, S.A.\ Molchanov,
Parabolic problems for the Anderson model II. Second-order asymptotics and structure of high peaks,
Probab.\ Theory Relat.\ Fields 111, 17--55, 1998.

\bibitem{G1999}
G.\ Grimmett,
\emph{Percolation} (2nd.\ ed.), 
Grundlehren der mathematischen Wissenschaften, 
Volume 321, Springer, Berlin, 1999.

\bibitem{K2016}
W.\ K\"onig,
\emph{The Parabolic Anderson Model},
Pathways in Mathematics, Birkh\"auser, 2016.

\bibitem{LP2016}
R.\ Lyons, Y.\ Peres,
\emph{Probability on Trees and Networks},
Cambridge Series in Statistical and Probabilistic Mathematics,
Cambridge University Press, New York, 2016.

\bibitem{BKS2018}
M.\ Biskup, W.\ K\"onig, R.S.\ dos Santos,
Mass concentration and aging in the parabolic Anderson model with doubly-exponential tails,
Probab.\ Theory Relat.\ Fields 171, 251--331, 2018. 

\bibitem{dhW2022}
F.\ den Hollander, D.\ Wang,
The parabolic Anderson model on a Galton-Watson tree revisited,
J. Stat. Phys. 189, no.1, Paper No. 8, 2022.

\bibitem{dhW2023}
F.\ den Hollander, D.\ Wang, 
The annealed parabolic Anderson model on a regular tree, 
Markov Proc. Relat. Fields (to appear).

\bibitem{AP2023}
E.\ Archer, A.\ Pein, Parabolic Anderson model on critical Galton–Watson trees in a Pareto environment, Stoch. Proc. Appl. 159, 34--100, 2023.

\bibitem{dHG1999}
J.\ G\"artner, F.\ den Hollander, Correlation structure of intermittency in the parabolic Anderson model, Probab. Theory Relat. Fields 114, 1--54, 1999.

\bibitem{JJ2002}
J.\ Jost, M.P.\ Joy, Spectral properties and synchronization in coupled map
lattices, Phys.Rev.E 65(1), 2002
\end{thebibliography}
\end{document}